\documentclass[graybox]{svmult}

\usepackage{type1cm}          
\usepackage{makeidx}          
\usepackage{graphicx}         
\usepackage{multicol}         
\usepackage[bottom]{footmisc} 
\usepackage{newtxtext}        %
\usepackage{newtxmath}        
\usepackage{indentfirst,csquotes}
\usepackage{listings}
\usepackage{mathtools}
\usepackage{adjustbox}
\usepackage{stmaryrd}
\usepackage{ulem}
\usepackage[most]{tcolorbox}
\usepackage{mdframed}
\usepackage{algorithm}
\usepackage{algpseudocode}
\usepackage{pbox}
\usepackage{xcolor}
\usepackage{hyperref}
\usepackage{etoolbox}
\usepackage{cleveref}

\input{preamble.sty}
\makeindex

\begin{document}

    \title*{Optimized Multilevel Sampling Methods \\ under Resource Constraints}
    \author{Niklas Baumgarten}
    \institute{Niklas Baumgarten \at Heidelberg University, \email{niklas.baumgarten@uni-heidelberg.de}}

    \maketitle

    \abstract{
    We present recent developments in multilevel
    sampling methods under resource constraints.
    Over the past 15 years, multilevel methods have
    become widely used for uncertainty quantification.
    However, scaling them to high-dimensional problems
    and high-performance computing (HPC) environments remains challenging.
    In this work, we discuss two algorithms designed to address these issues:
    the budgeted Multilevel Monte Carlo (MLMC) method
    and the Multilevel Stochastic Gradient Descent (MLSGD) method.
    We demonstrate their effectiveness on HPC systems
    under consideration of the available computational resources
    for applications in forward uncertainty quantification (UQ)
    and optimal control (OC) under uncertainty.
}

    \section{Introduction}\label{sec:introduction}

Collecting and merging large amounts of data in complex
computational systems present significant challenges,
especially when there is a limited amount of memory and processing power.
For tasks like uncertainty quantification (UQ) and optimal control (OC)
of partial differential equations (PDEs),
these computational limits restrict the precision
of the achievable results using numerical methods.

We present two strongly related algorithms that can operate
within the imposed constraints by the computer system and
achieve optimal results in terms of memory and CPU-time usage.
The first algorithm, the budgeted Multilevel Monte Carlo (MLMC) method~\cite{baumgarten2025budgeted},
based on the initial work
of~\cite{giles2008multilevel, cliffe2011multilevel, barth2011multi, giles2015multilevel},
is designed to compute full field estimates of PDE systems under uncertainty.
Full field estimates, e.g.~of the mean field, can find application
in weather forecasting, climate simulations, but also in material
science to comprehensively describe material properties
or in other related fields in engineering and physics.

The second algorithm, the Multilevel Stochastic Gradient Descent (MLSGD)
method~\cite{baumgarten2025multilevel},
finds the optimal control to a PDE system under uncertainty.
Its development is inspired by the work of
\cite{geiersbach2020stochastic, guth2023multilevel, van2019robust}
and uses the full field MLMC method from~\cite{baumgarten2025budgeted}
to estimate the gradients within the optimization process.

Both algorithms achieve exceptional performance by
leveraging a multiindex finite element (FE) mesh on distributed memory,
extending the work of~\cite{baumgarten2023fully, baumgarten2024fully},
and by applying currying techniques to approximate
solutions of multiple coupled PDEs.
In the context of full field estimation,
these coupled PDEs involve the sampling of Gaussian Random Fields
(GRFs) with the techniques of~\cite{lindgren2011explicit, kutri2024dirichlet},
solving a subsurface diffusion problem,
and then computing mass transport in the resulting flux field.
The coupling arises because each stage's output defines the input for the next.
For the optimal control problem,
the PDE coupling similarly begins with GRF sampling,
followed by solving the state and adjoint equations to compute gradients,
as outlined in~\cite{geiersbach2019projected}.

The algorithms are implemented in the latest version~\cite{wieners2025mpp350}
of the software M++
(see~\cite{baumgarten2021parallel} for an introduction of the software)
and developed following a rigorous process detailed
in~\cite{baumgarten2025continuous}.
M++ plays a central role in the HPC project associated
with this report and is applied across various scientific domains.
In geoscience, it supports full waveform inversion for
seismic imaging~\cite{bohlen2020visco} and gas dynamics simulations for
carbon capturing~\cite{knodel2022global}.
In materials science, it is used to model nonlinear
solid mechanics~\cite{bayat2018numerical}
and dislocation dynamics~\cite{schulz2019mesoscale, wagner2019discontinuous},
and M++ is also applied in life sciences to simulate cardiovascular
processes~\cite{frohlich2023numerical, froehlich2022cardiac, gerach2021electro}.

Here, we focus on applications of M++ in UQ and OC,
where our core approach involves solving knapsack problems.
These problems naturally arise when we try to minimize errors,
such as the root mean square error (RMSE),
under a fixed computational budget given by the available CPU-time
$\rP \cdot \rT_0$ and memory $\mathrm{Mem}_0$.
The main challenge lies in allocating limited and distributed resources
across multiple resolution levels, samples, and coupled PDEs.
To address this, we employ distributed dynamic programming techniques,
detailed in~\cite{baumgarten2024fully, baumgarten2025budgeted, baumgarten2025multilevel},
and demonstrate their performance in our numerical experiments.
The report introduces the forward UQ and OC problems in~\ref{sec:problem-statement},
and outlines the methodology and numerical results in
Section~\ref{sec:numerical-experiments}.
Section~\ref{sec:outlook} presents an outlook on future developments.

    \section{Problem Statements} \label{sec:problem-statement}

This report addresses two classes of problems:
forward UQ for PDEs and PDE-constrained OC problems under uncertainty.
Both are posed on a bounded polygonal domain $\cD \subset \RR^d$, $d \in \tset{2, 3}$,
and potentially also over the time interval $[0, T]$.
The uncertainty is modeled using measures on the complete
probability space $(\Omega, \mathcal{F}, \PP)$,
and the solutions are sought in appropriate Bochner spaces.
These are composed of the probability space and the function spaces
$V \subset \rL^2(\cD)$ and $W \subset \rL^2(\cD)$
over the spatial domain and possibly also over the time interval
(see~\cite{baumgarten2025budgeted, baumgarten2025multilevel}
for definitions, analysis and further details).

\begin{problem}[Forward Uncertainty Quantificiaton]
    \label{problem:forward-uq}
    Goal is to estimate statistical moments, e.g.~mean fields $\EE[\bu] \in V$,
    of the solution $\bu \colon \Omega \times \cD \times [0, T] \rightarrow \RR$
    to a well-posed PDE system with randomly distributed input data, i.e.,
    systems of the form
    \begin{equation}
        \label{eq:generic-pde-model}
        \cG[\omega] \bu(\omega, \bx, t) = \bz(\bx, t)
    \end{equation}
    where $\omega \in \Omega$ corresponds to an outcome
    in the probability space, $\bx \in \cD$ denotes the spatial location
    and $t \in [0, T]$ represents a point in time.
    The right-hand side $\bz(\bx, t)$ is a forcing term to the PDE system,
    and the operator $\cG[\omega]$ is a generic (potentially non-linear) operator
    influenced by random events $\omega \in \Omega$,
    representing the coupled PDE models.
\end{problem}

\begin{figure}
    \centering

    \includegraphics[width=0.215\textwidth]{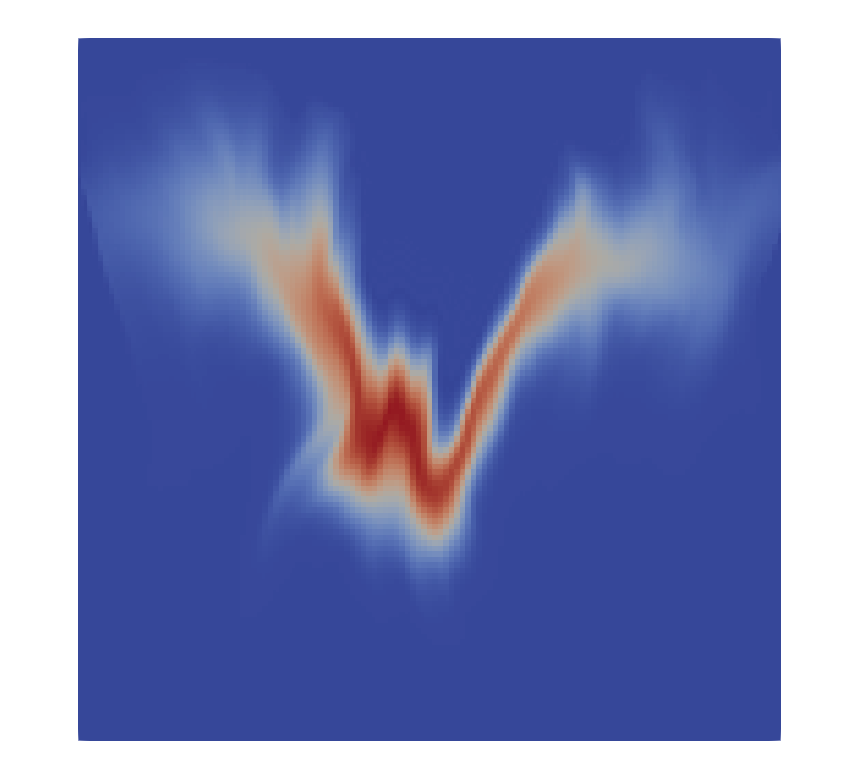}
    \hspace{-5mm}
    \includegraphics[width=0.215\textwidth]{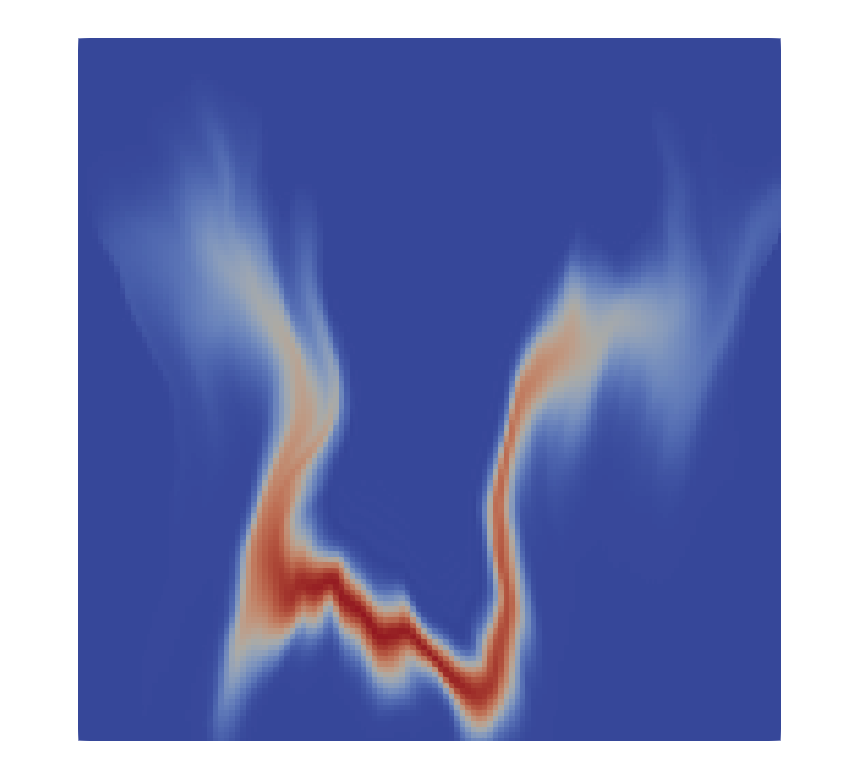}
    \hspace{-5mm}
    \includegraphics[width=0.215\textwidth]{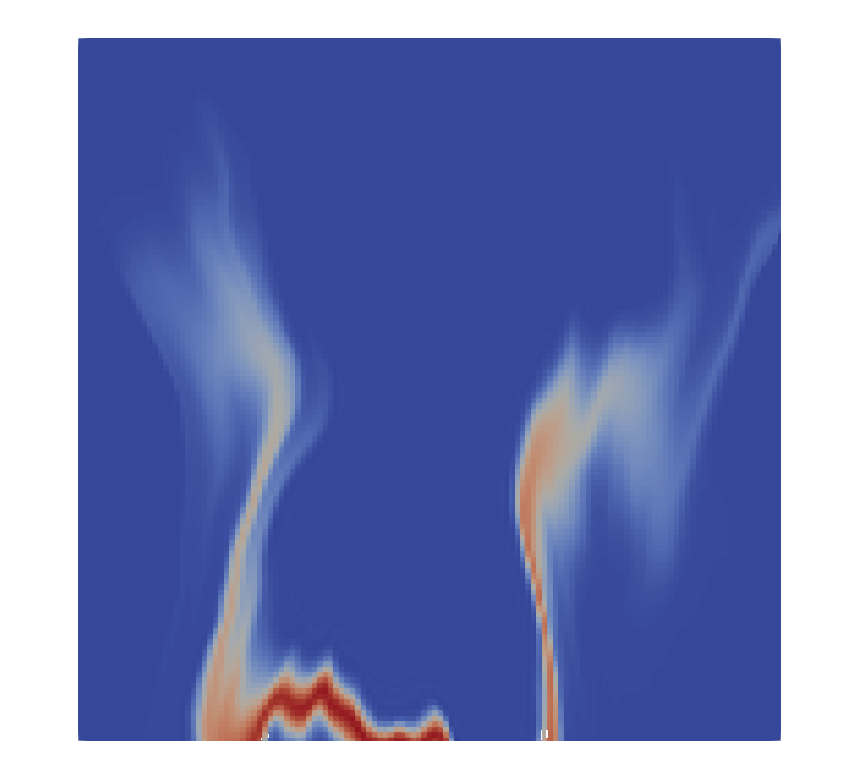}
    \hspace{-5mm}
    \includegraphics[width=0.215\textwidth]{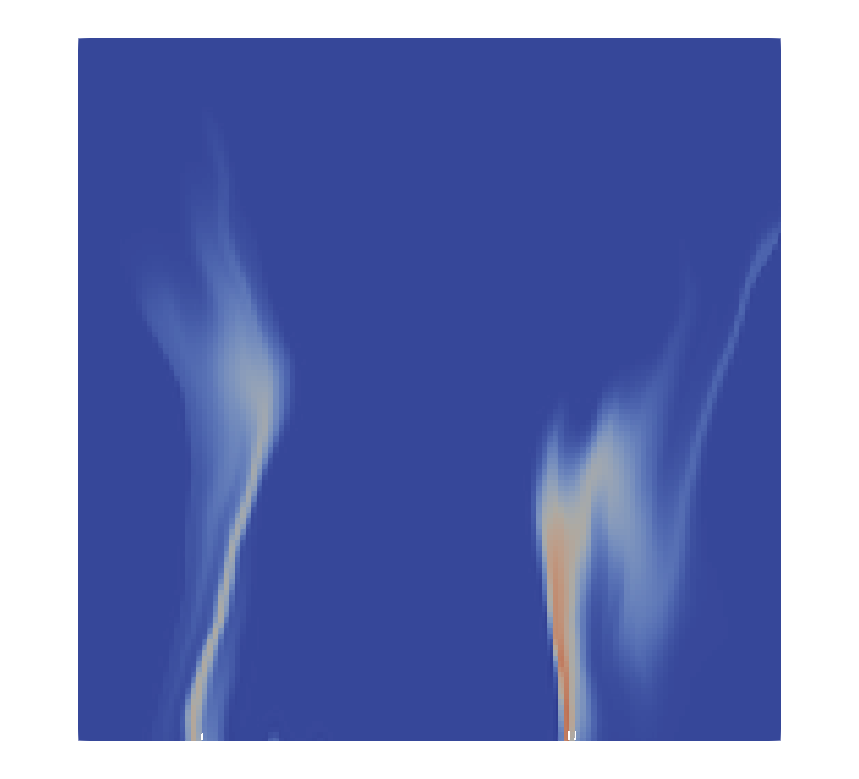}
    \hspace{-5mm}
    \includegraphics[width=0.215\textwidth]{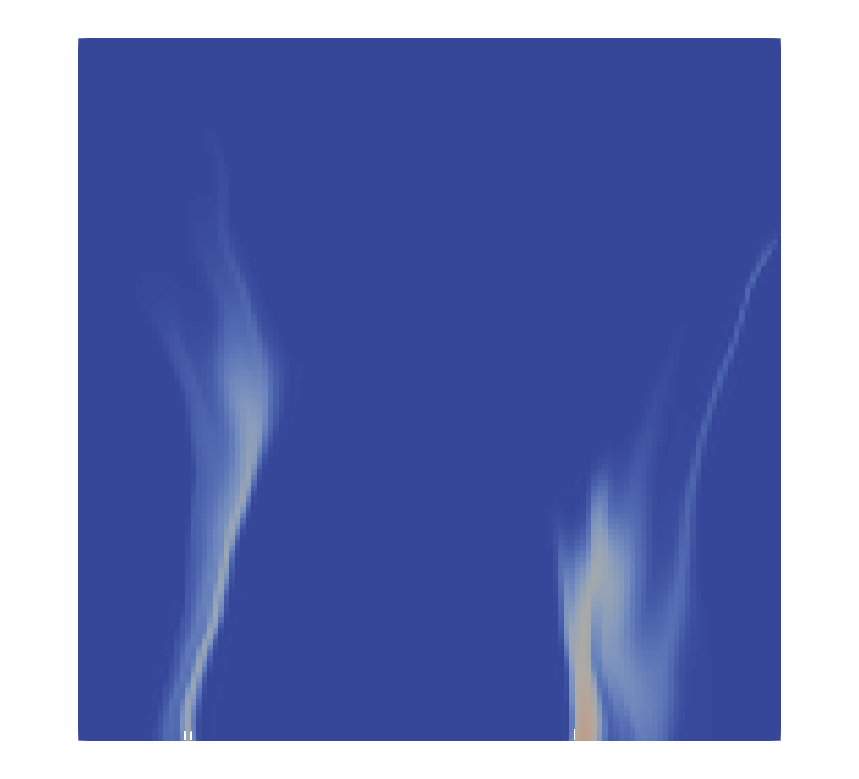}
    \vspace{-3mm}
    \caption{Temporal development of one sample following~\eqref{eq:hyperbolic-transport}.}
    \label{fig:sample1}

    \includegraphics[width=0.215\textwidth]{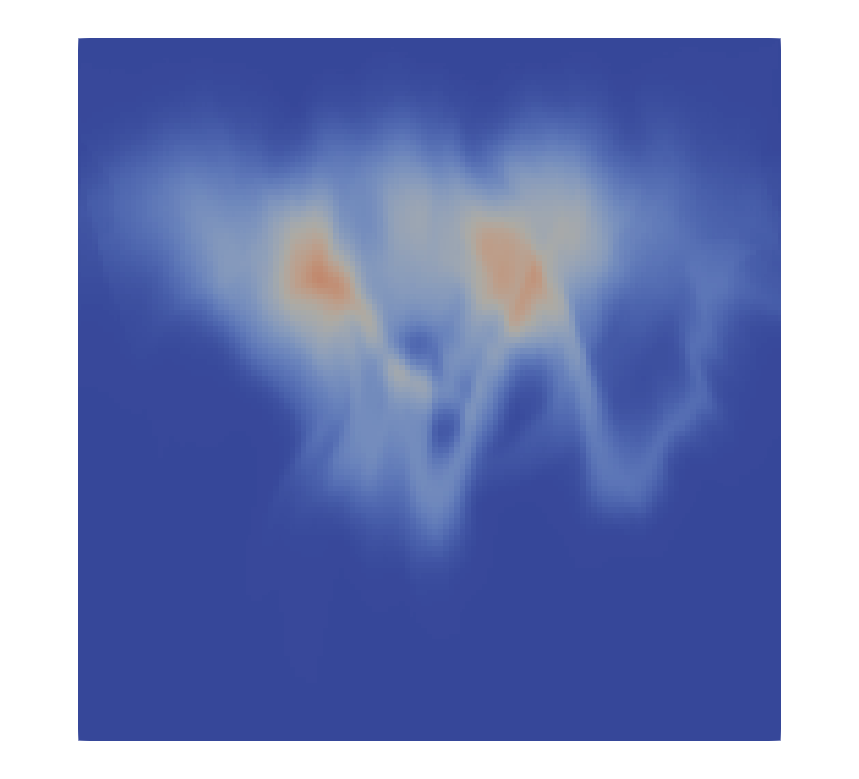}
    \hspace{-5mm}
    \includegraphics[width=0.215\textwidth]{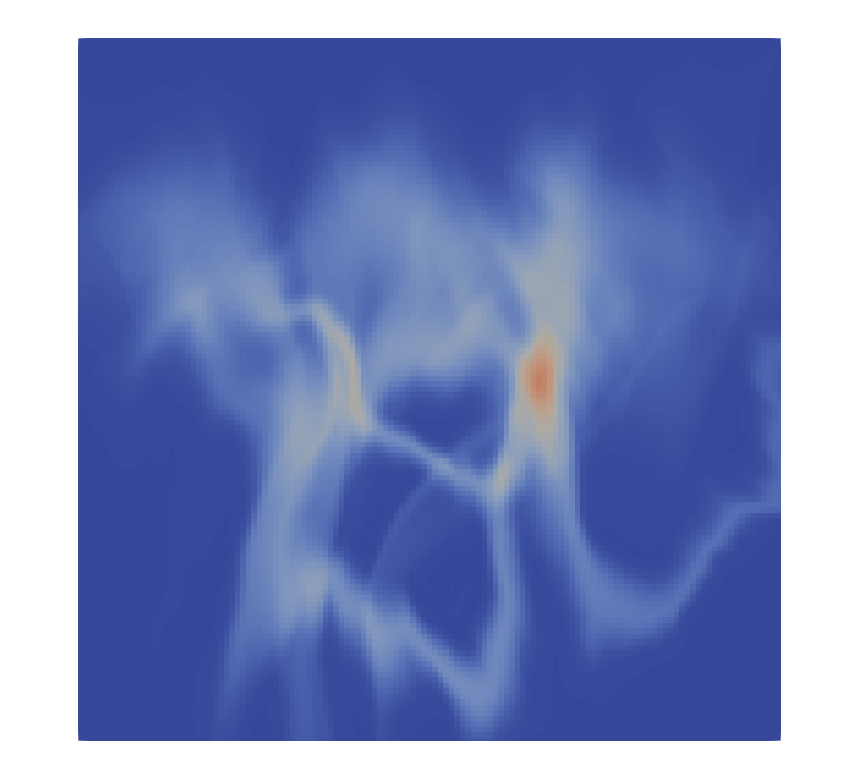}
    \hspace{-5mm}
    \includegraphics[width=0.215\textwidth]{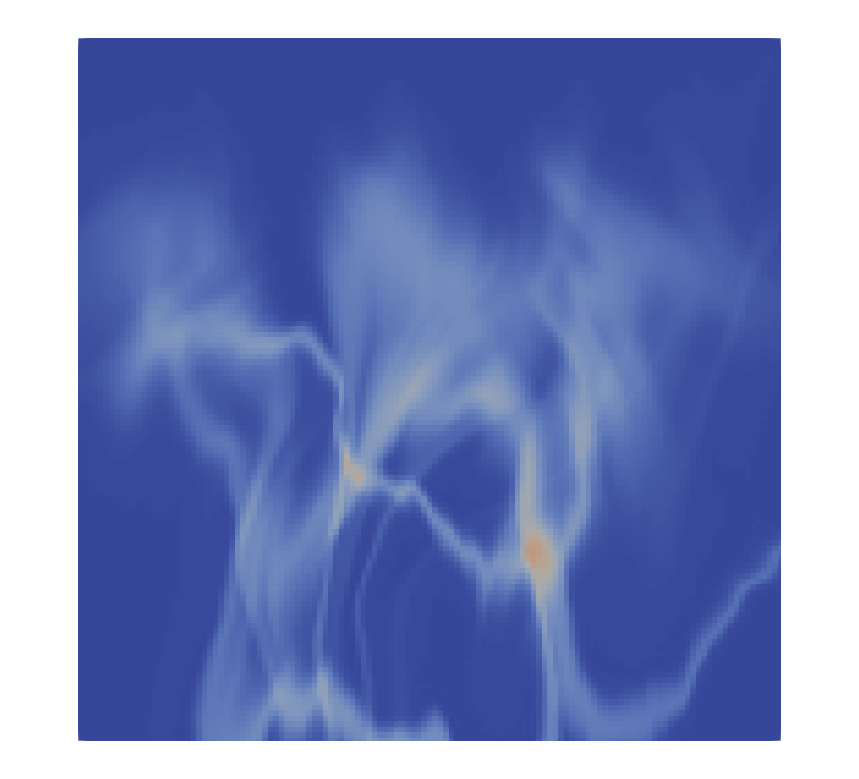}
    \hspace{-5mm}
    \includegraphics[width=0.215\textwidth]{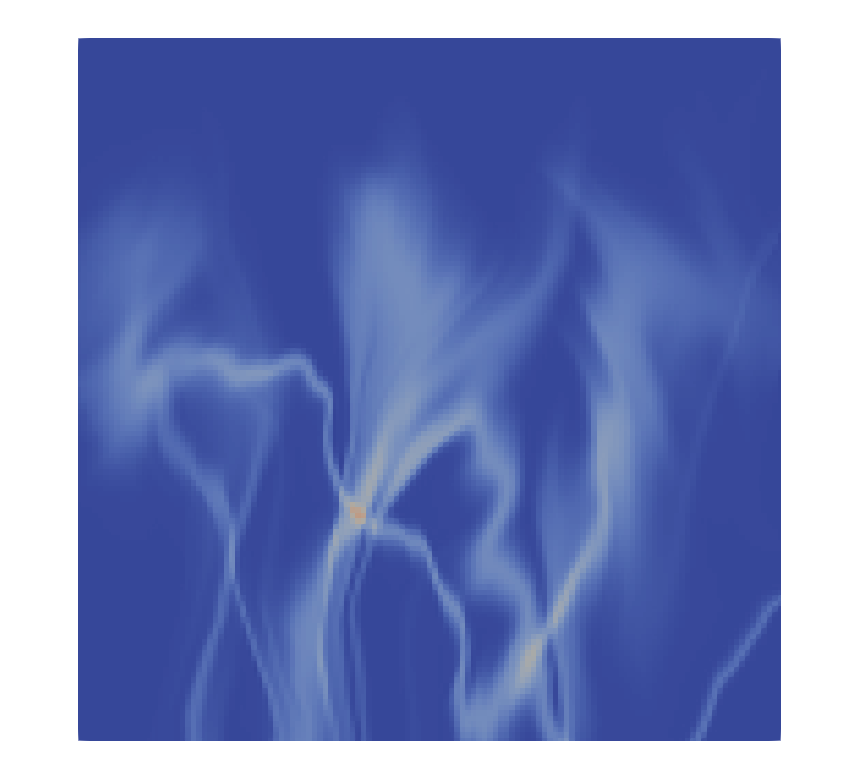}
    \hspace{-5mm}
    \includegraphics[width=0.215\textwidth]{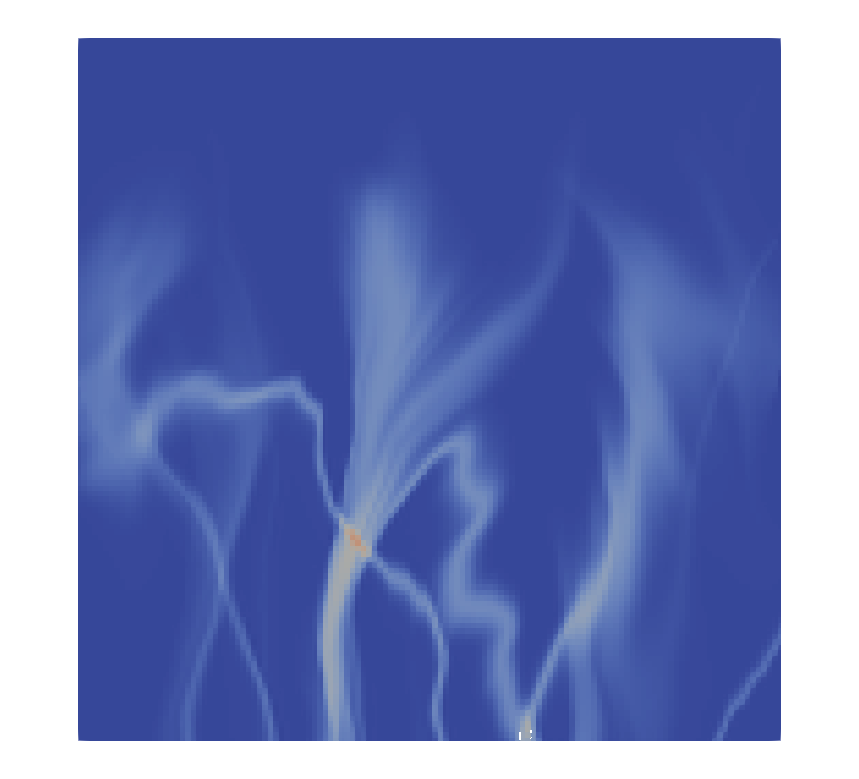}
    \vspace{-3mm}
    \caption{Four samples averaged of the transport system~\eqref{eq:hyperbolic-transport}.}
    \label{fig:mean-para}

    \includegraphics[width=0.215\textwidth]{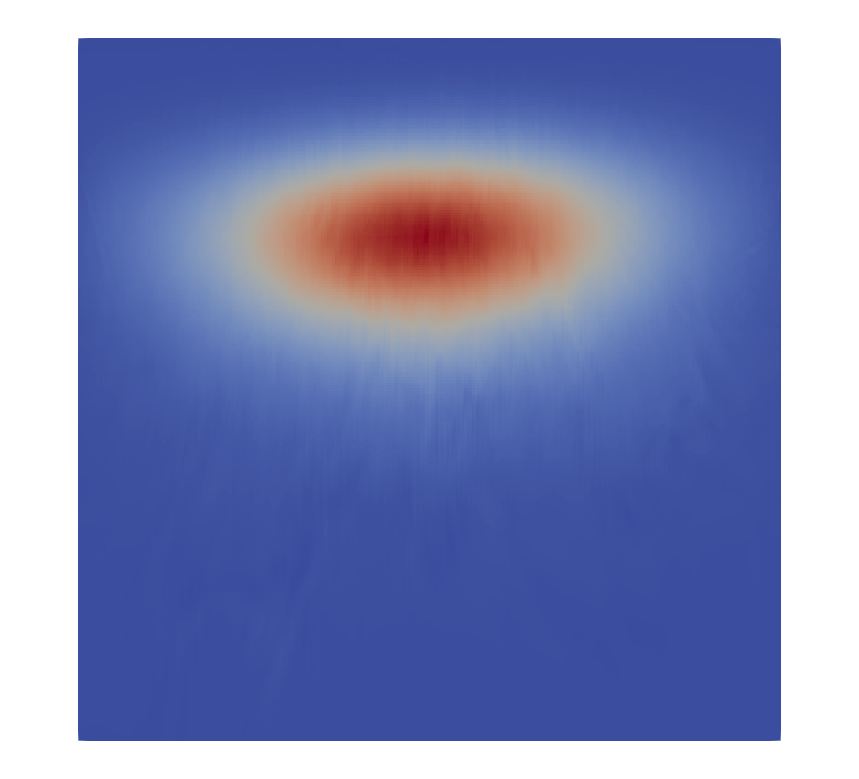}
    \hspace{-5mm}
    \includegraphics[width=0.215\textwidth]{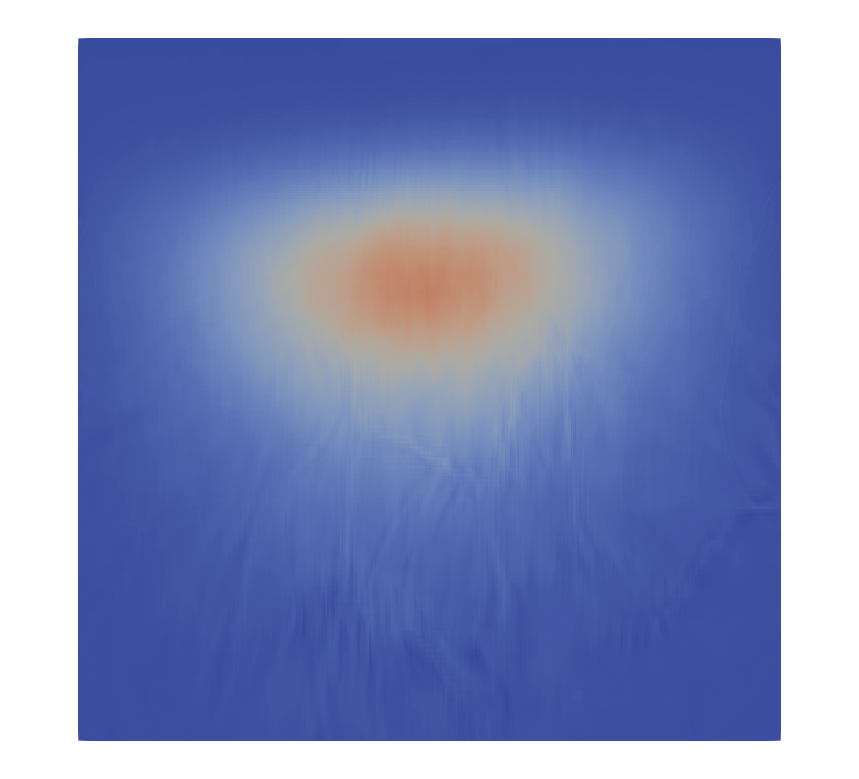}
    \hspace{-5mm}
    \includegraphics[width=0.215\textwidth]{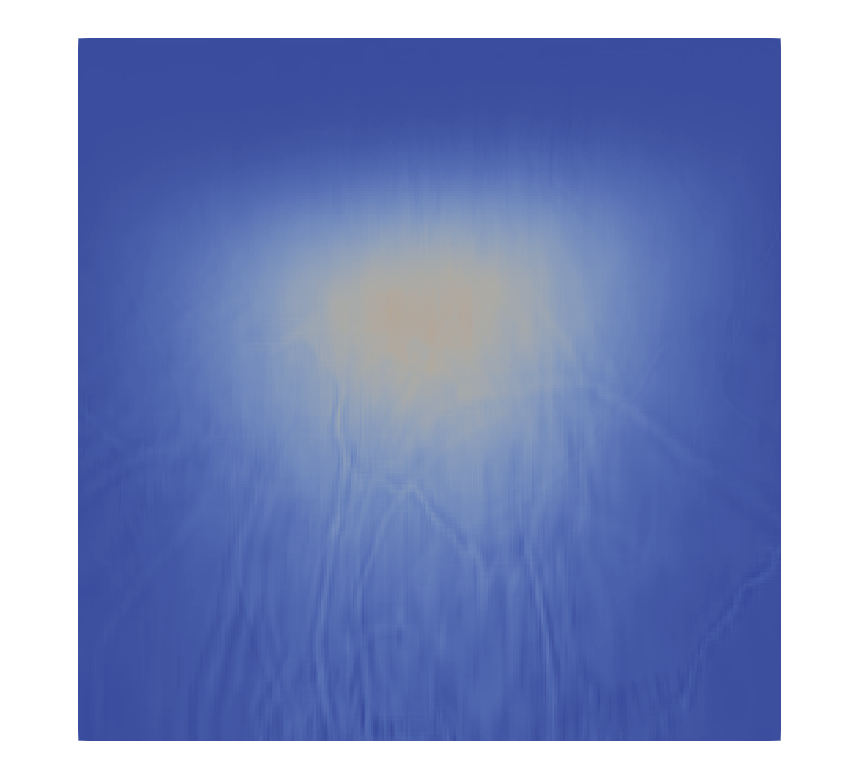}
    \hspace{-5mm}
    \includegraphics[width=0.215\textwidth]{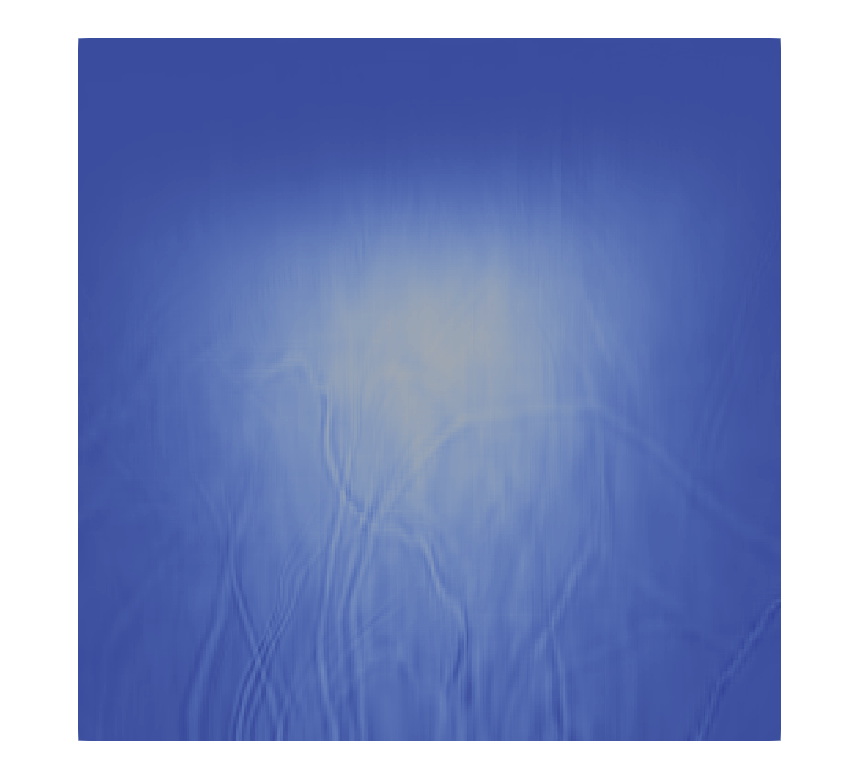}
    \hspace{-5mm}
    \includegraphics[width=0.215\textwidth]{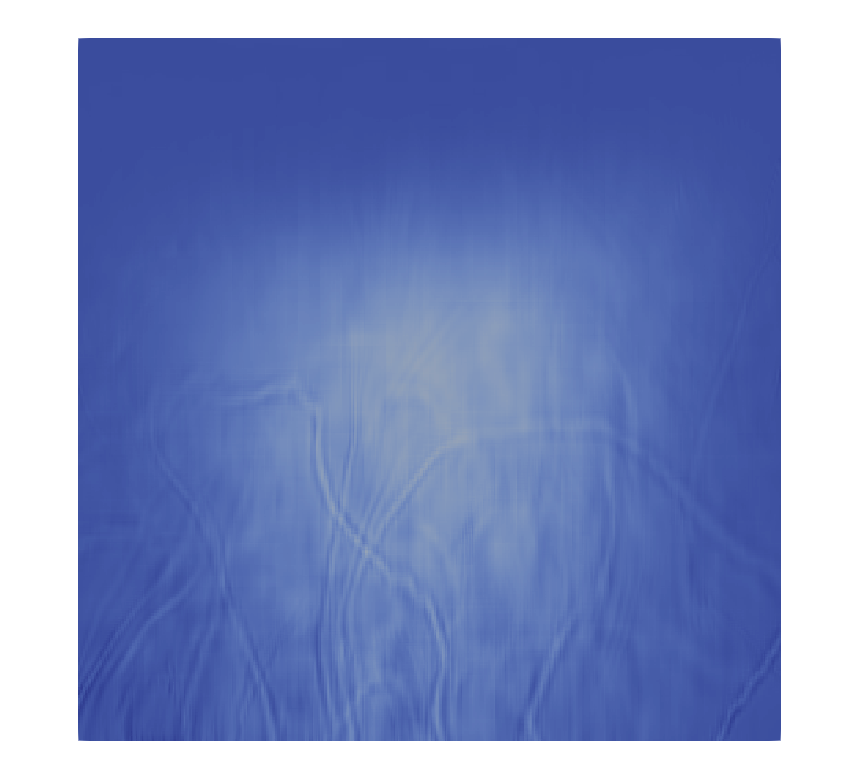}
    \vspace{-3mm}
    \caption{Multilevel estimate of the mean field $\EE[\bu]$ with estimator~\eqref{eq:mlmc-estimator}.}
    \label{fig:mean-intermediate}

    \includegraphics[width=0.215\textwidth]{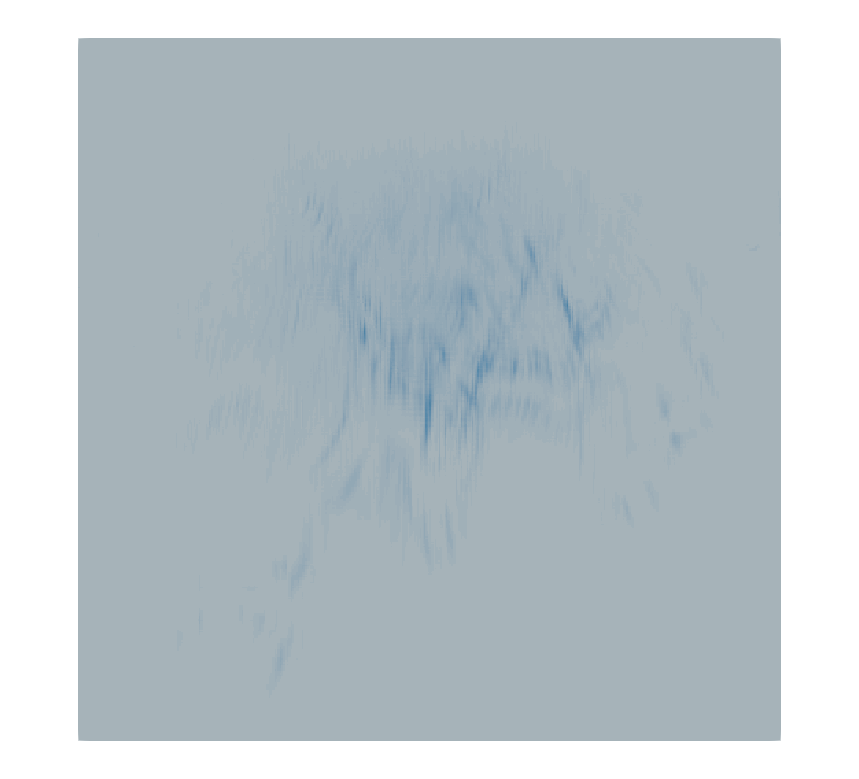}
    \hspace{-5mm}
    \includegraphics[width=0.215\textwidth]{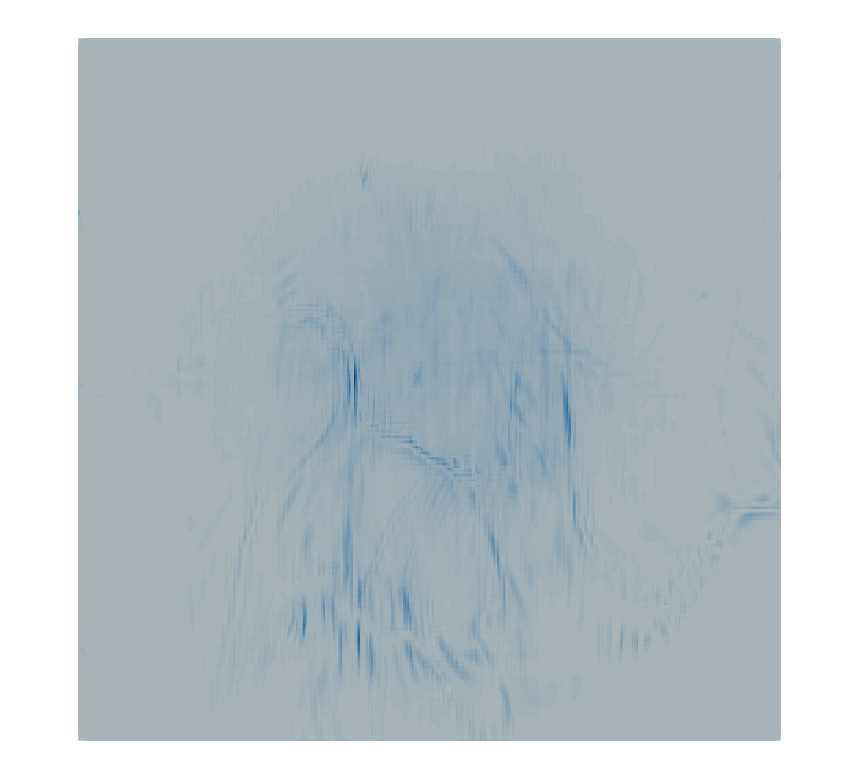}
    \hspace{-5mm}
    \includegraphics[width=0.215\textwidth]{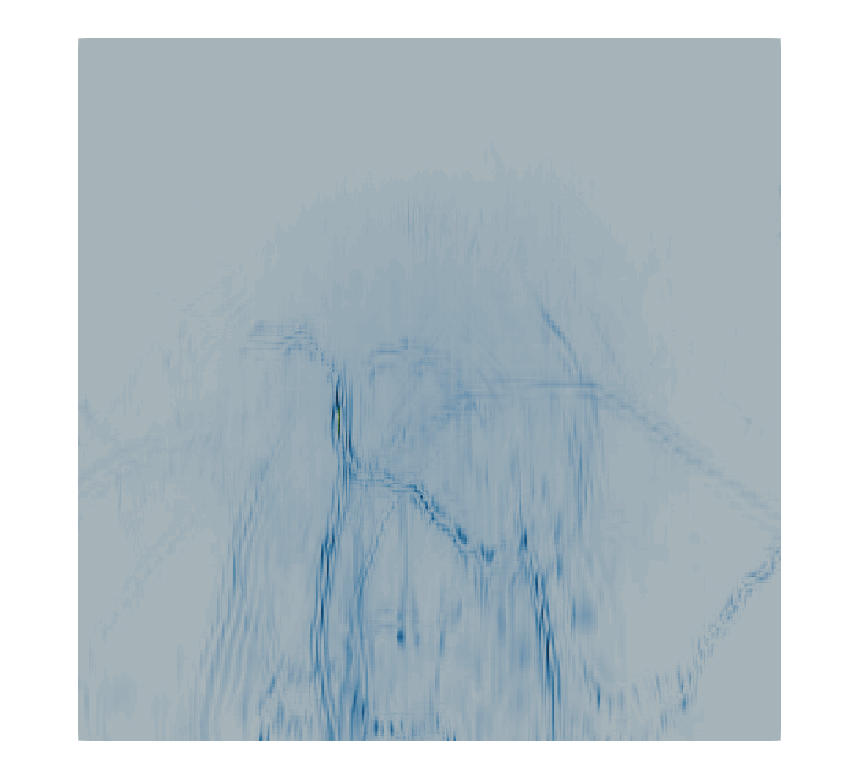}
    \hspace{-5mm}
    \includegraphics[width=0.215\textwidth]{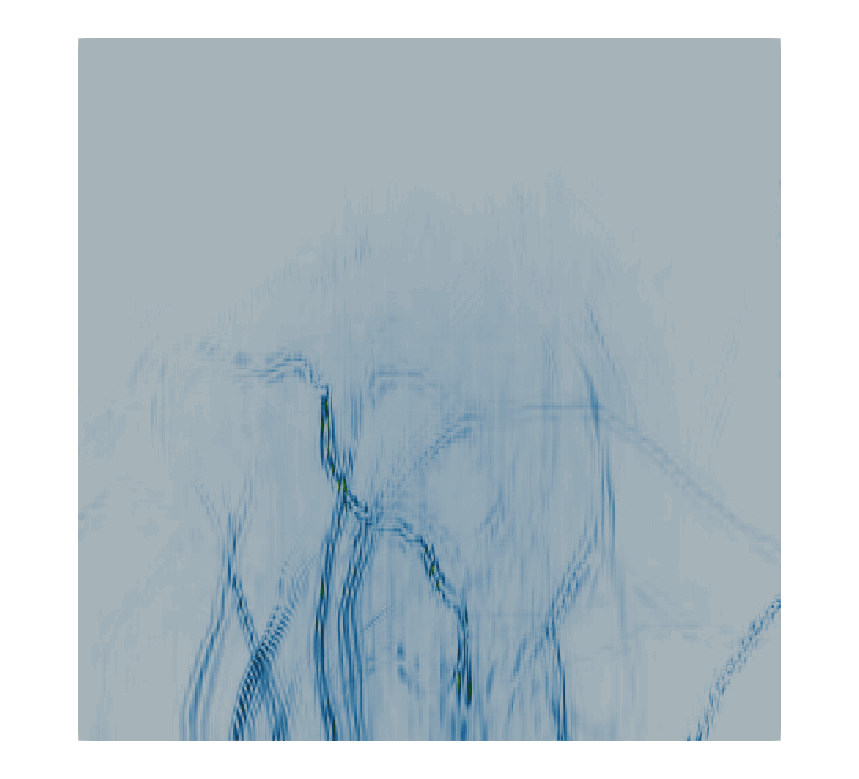}
    \hspace{-5mm}
    \includegraphics[width=0.215\textwidth]{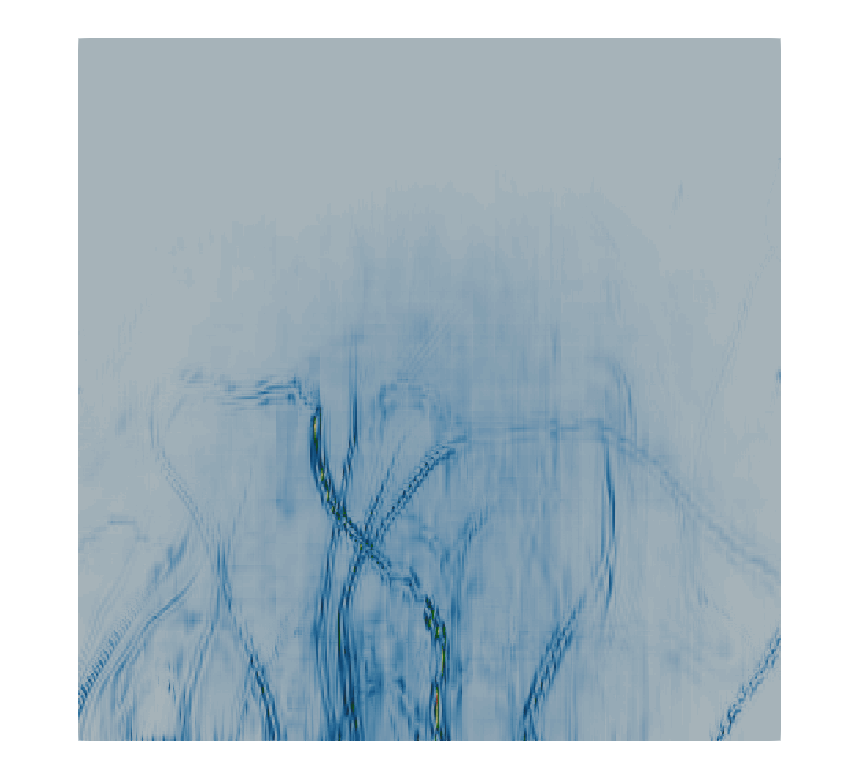}
    \vspace{-3mm}
    \caption{Multilevel marginal sample variance estimate over time after completion of the estimation.}
    \label{fig:svar-intermediate}
\end{figure}

As a concrete application of Problem~\ref{problem:forward-uq},
we consider the computation of the mass distribution
$\bu \colon \Omega \times \cD \times [0, T] \rightarrow \RR$
governed by the hyperbolic transport system:
\begin{equation}
    \label{eq:hyperbolic-transport}
    \pdeProblem{
        \partial_t \bu(\omega, \bx, t) + \div(\bq(\omega, \bx) \bu(\omega, \bx, t)) &=& 0 & \text{on } \cD \times (0, T] \\
        \bu(\omega, \bx, t) &=& \bu(\bx, 0) & \text{on } \cD
    }
\end{equation}
Here, $\bq(\omega, \bx)$ is the flux field derived from
a subsurface diffusion problem with log-normal coefficients,
based on GRF realizations obtained via SPDE
sampling~\cite{lindgren2011explicit, kutri2024dirichlet}.
The operator $\cG[\omega]$ thus is the composition of GRF sampling,
subsurface diffusion, and the hyperbolic transport PDE system.
For further details, see~\cite{baumgarten2025budgeted};
Figures~\ref{fig:sample1}-\ref{fig:svar-intermediate} illustrate the computational task.

Figure~\ref{fig:sample1} shows the temporal evolution
of a single sample $\bu_\ell^{(m)} \in V_\ell$, where $m$ denotes the sample index,
$\ell$ denotes the FE level and $V_\ell$ the associated FE space.
The randomly drawn input flux field $\bq^{(m)}_\ell$
clearly distorts the initial spatial
Gaussian mass distribution $\bu(\bx, 0)$ over time.
The second Figure~\ref{fig:mean-para}, presents the average of four samples
$\tfrac{1}{4}\sum_{m=1}^4 \bu^{(m)}_\ell$,
which does not yet exhibit the desired properties
of a robust full field estimate but serves as
an intermediate result which can be computed in parallel.
The third Figure~\ref{fig:mean-intermediate},
shows a multilevel estimate $E^{\text{ML}}[\bu_L]$
of the mean field $\EE[\bu]$,
\begin{equation}
    E^{\text{ML}}[\bu_L] \coloneqq E^{\text{MC}}_{M_0}[\bu_0]
    + \sum_{\ell=1}^L E^{\text{MC}}_{M_\ell}[\bu_\ell - \bu_{\ell - 1}]
    \approx \EE[\bu]
    \label{eq:mlmc-estimator}
\end{equation}
where all $\bu_\ell$ are FE approximations of the solution and $E^{\text{MC}}_{M_\ell}$
are standard Monte Carlo (MC) estimates,
i.e.~$E^{\text{MC}}_{M_0}[\bu_\ell] \coloneqq \tfrac{1}{M_\ell} \sum_{m=1}^{M_\ell} \bu_\ell^{(m)}$,
with $M_\ell$ independent samples on different mesh levels $\ell$.

The telescoping sum~\eqref{eq:mlmc-estimator}
enables the efficient estimation of quantities through
variance reduction, which will lead as we will
see in Section~\ref{sec:numerical-experiments}, to drastic
speedups compared to a standard single-level MC estimation.
The final row, Figure~\ref{fig:svar-intermediate},
illustrates the marginal sample variance field
after completing the multilevel estimation.
Thus, this plot indicates how large the error is at
each point in space and time.
We further remark, that we use discontinuous Galerkin (dG)
methods with upwind fluxes
and use the implicit midpoint rule for time discretization to compute the samples $\bu_\ell^{(m)}$.

We now turn to the second problem class:
the PDE-constrained optimal control problem.
The goal is to find an optimal control $\bz \in W$
such that the corresponding state solution of the PDE is as
close as possible to a given target $\bd \in W$ in expectation.

\begin{problem}[Optimal Control under Uncertainty]
    \label{problem:ocp}
    Given the desired target state $\bd \in W$ and a cost factor $\lambda\geq 0$, the goal is to
    find the optimal control $\bz \in Z$, such that
    \begin{equation}
        \label{eq:ocp-objective}
        \min_{\bz \in Z} J(\bz) \coloneqq \EE \squarelr{j(\cdot, \bz)} \quad \text{with} \quad
        j(\omega, \bz) \coloneqq
        \tfrac{1}{2} \norm{\bu[\omega] - \bd}^2_{W}
        + \tfrac{\lambda}{2} \norm{\bz}_{W}^2
    \end{equation}
    under the constraint that the state $\bu \in \rL^2(\Omega, V)$ is the solution of
    \begin{equation}
        \label{eq:ocp-constraint}
        \cG[\omega] \, \bu(\omega, \bx) = \bz(\bx)
    \end{equation}
    with $\cG[\omega]$ representing again the uncertain PDE system.
\end{problem}

\begin{figure}
    \begin{center}
        \raisebox{0.55in}{\rotatebox{90}{$h_\ell=2^{-7}$}}
        \includegraphics[trim={7cm 1.8cm 7cm 0}, clip, width=0.3\textwidth]{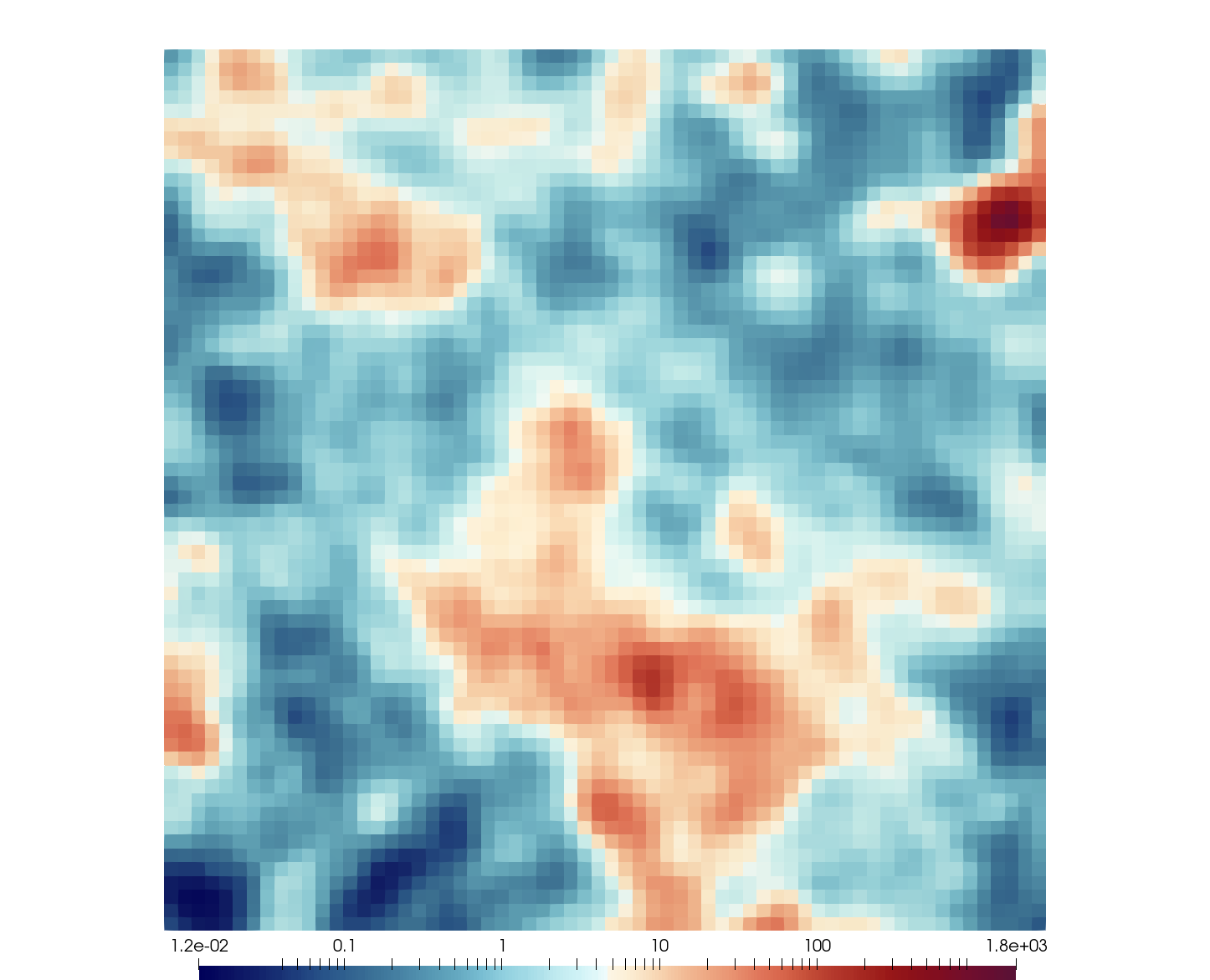}
        \includegraphics[trim={7cm 1.8cm 7cm 0}, clip, width=0.3\textwidth]{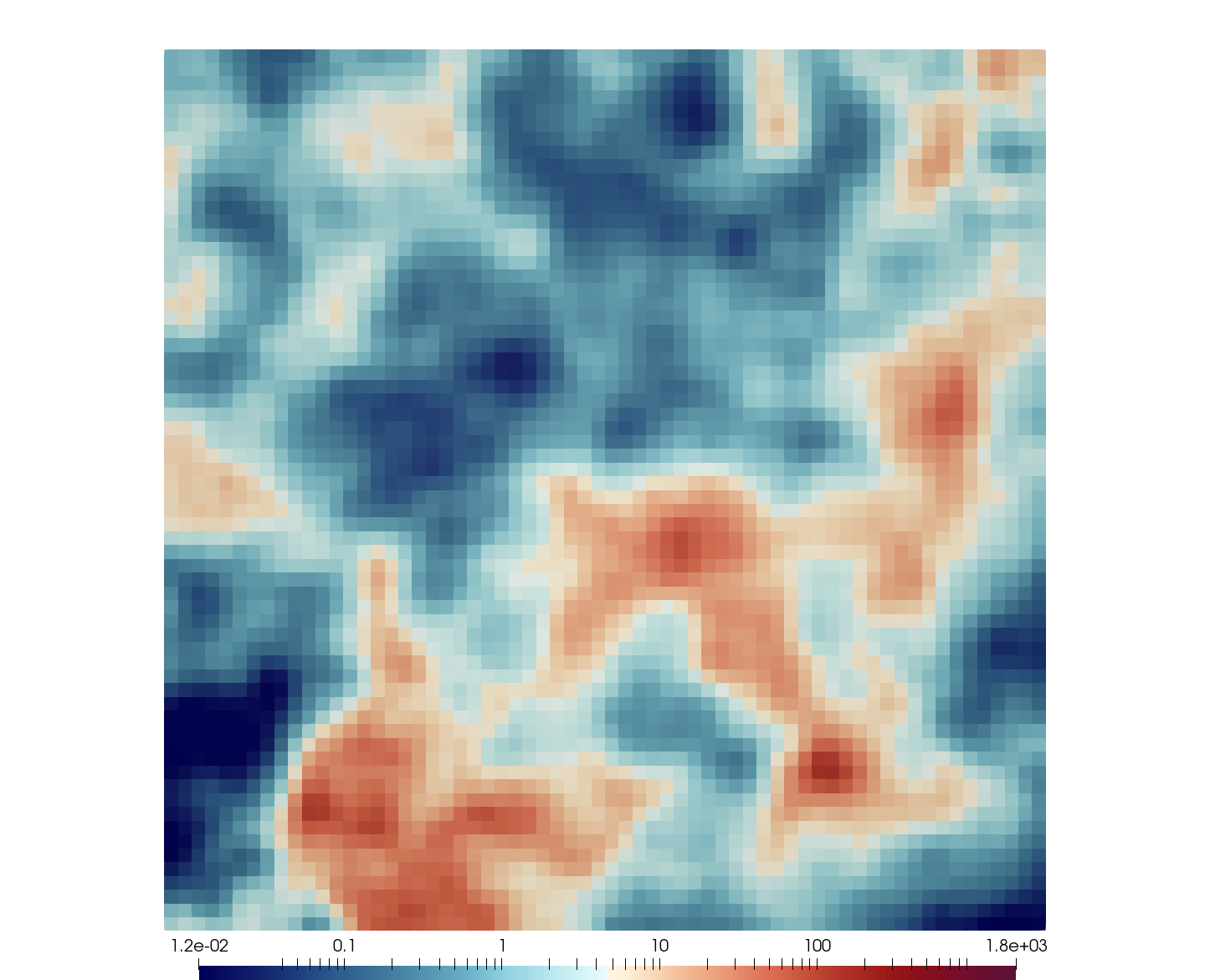}
        \hspace{0.1mm}
        \includegraphics[trim={7cm 1.8cm 7cm 0}, clip, width=0.3\textwidth]{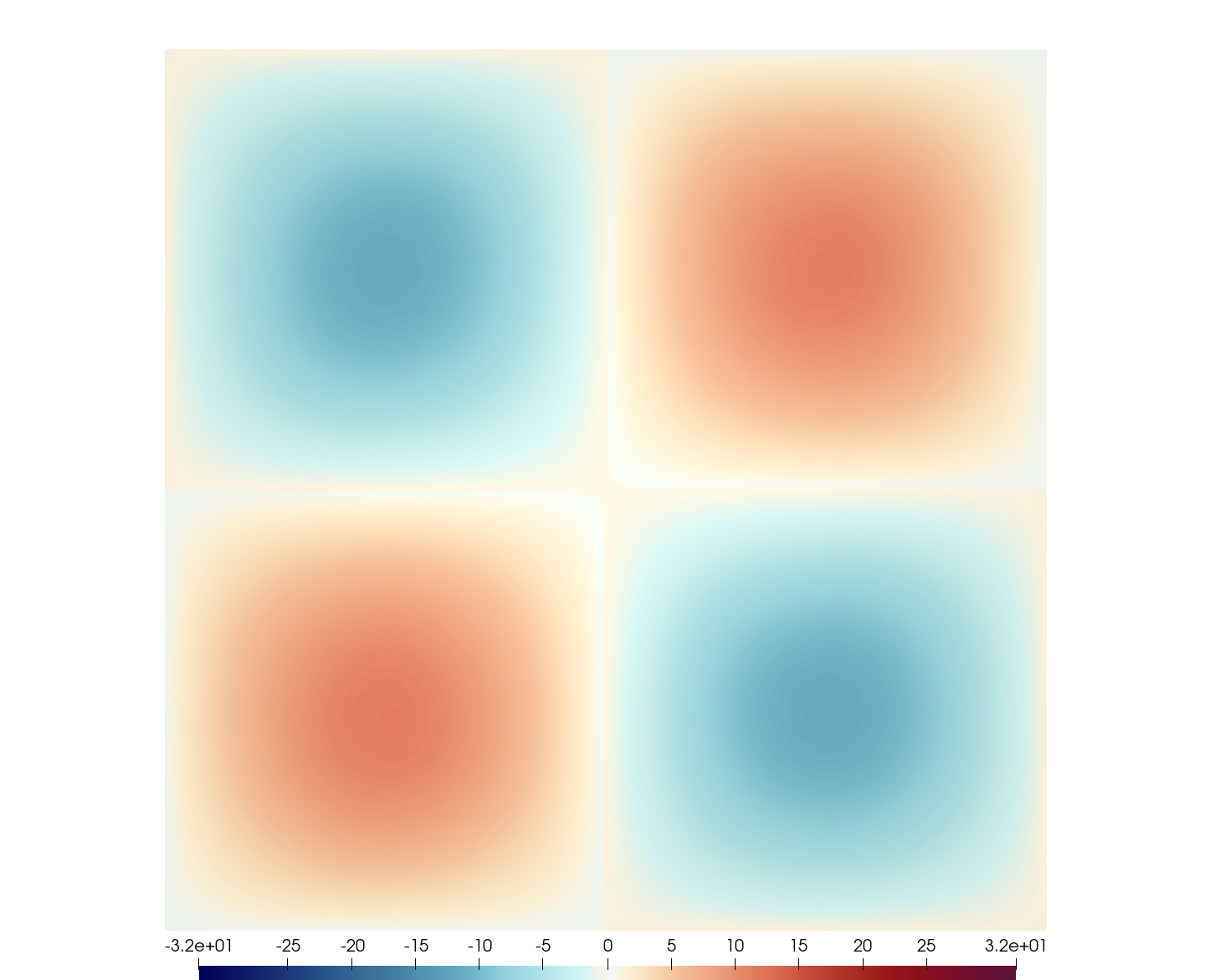}
        \raisebox{1.0in}{\rotatebox{-90}{$k=10$}}

        \vspace*{-1.5mm}

        \raisebox{0.55in}{\rotatebox{90}{$h_\ell=2^{-8}$}}
        \includegraphics[trim={7cm 2cm 7cm 0.4cm}, clip, width=0.3\textwidth]{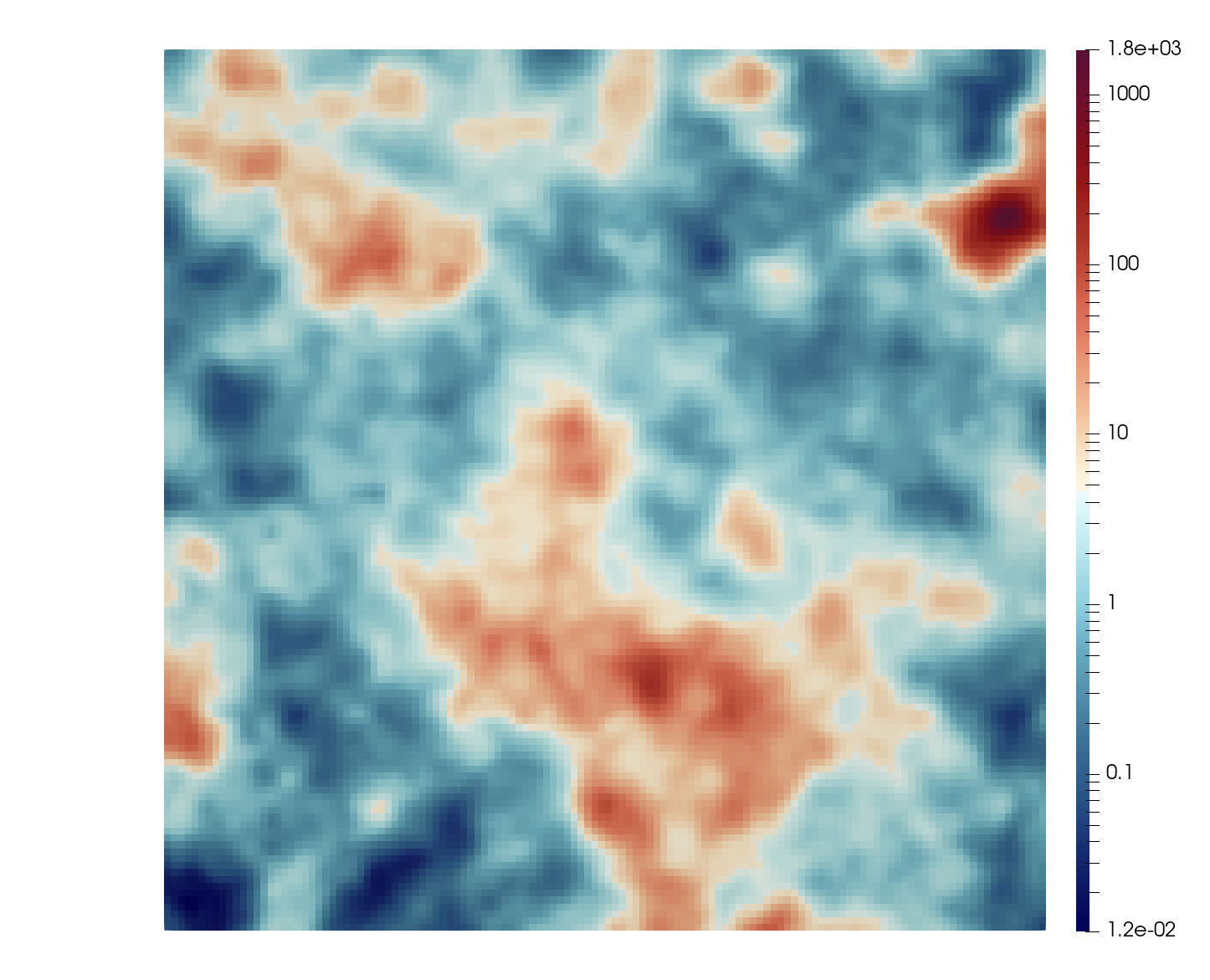}
        \includegraphics[trim={7cm 2cm 7cm 0.4cm}, clip, width=0.3\textwidth]{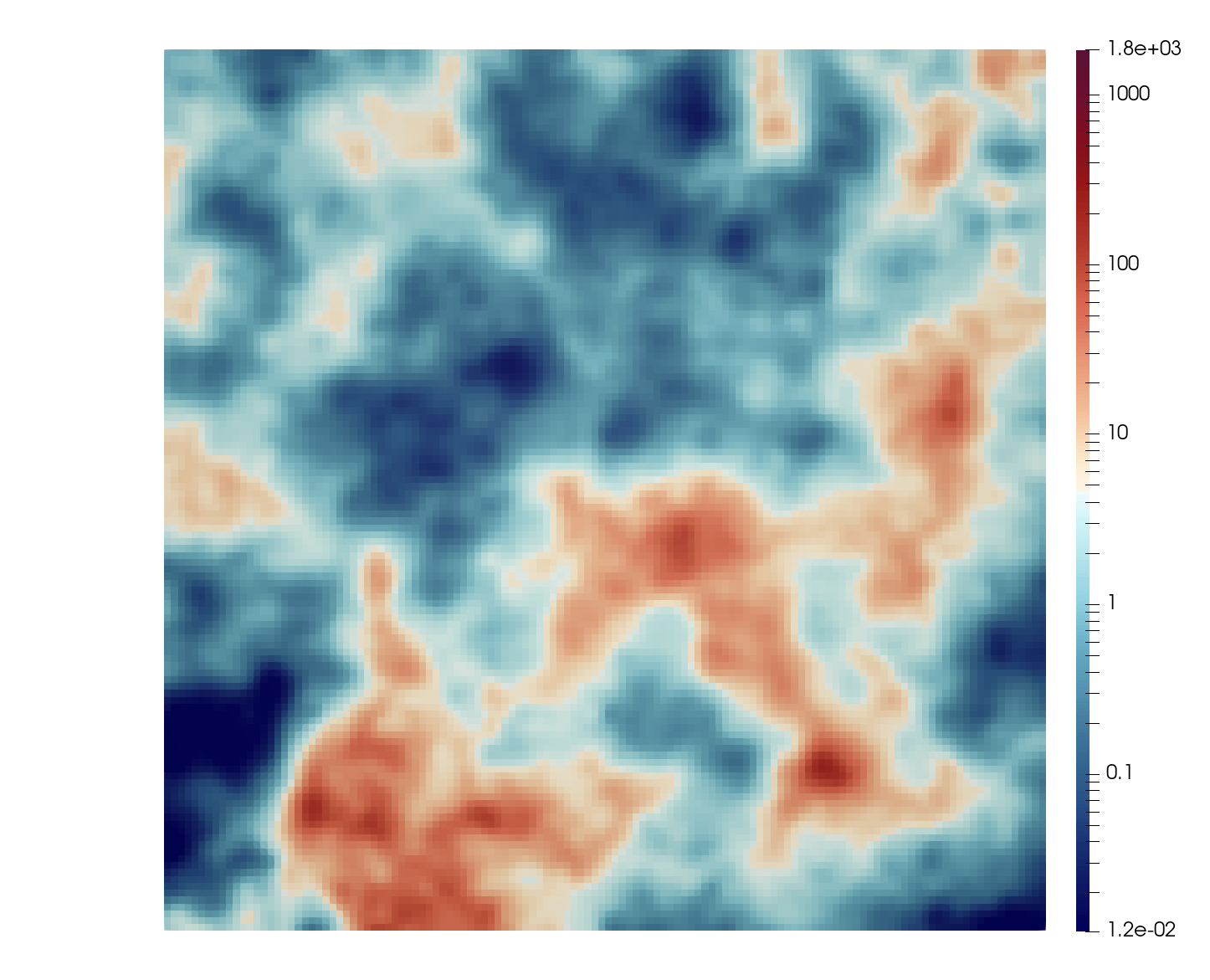}
        \hspace{0.1mm}
        \includegraphics[trim={7cm 2cm 7cm 0.4cm}, clip, width=0.3\textwidth]{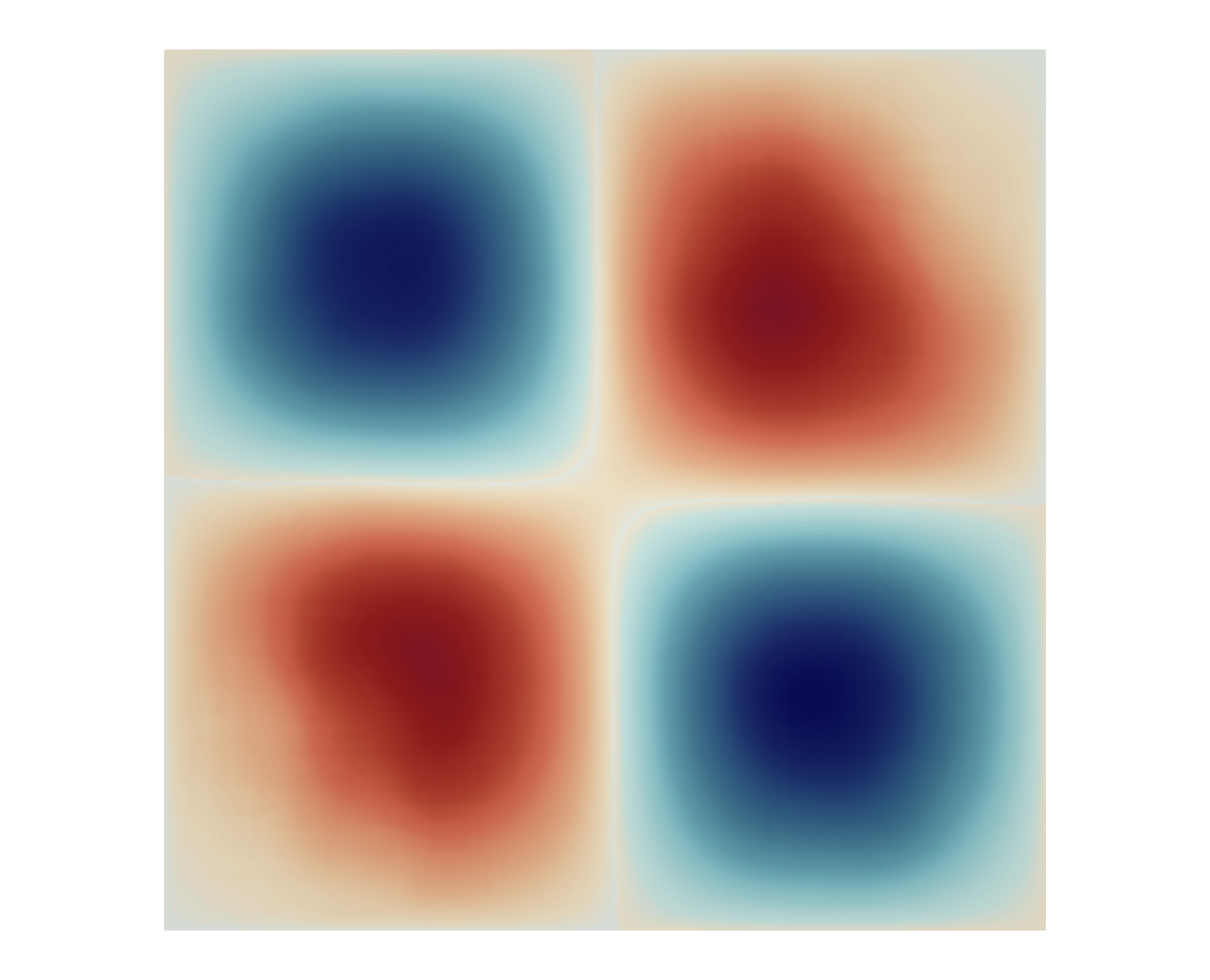}
        \raisebox{1.0in}{\rotatebox{-90}{$k=100$}}
    \end{center}
    \vspace*{-0.4cm}
    \caption{Left to right: Two GRF samples $\by_\ell^{(m)}$ on different mesh diameters;
    Computed control $\bz^{(k)}_\ell$ after $k=10$ and $k=100$ iterations of~\eqref{eq:bsgd-iteration}.}
    \label{fig:log-normal-fields-and-control}
\end{figure}

In contrast to the previous problem, the PDE system
is not time-dependent anymore and we consider for
$\cG[\omega]$ an eliptic diffusion equation on
$W=\rL^2(\cD)$ and $V=\rH^1_0(\cD)$ with
log-normal coefficients
and homogeneous Dirichlet boundary conditions:
\begin{equation}
    \pdeProblem{
        -\div\big(\exp(\by(\omega, \bx)) \nabla \bu(\omega, \bx) \,\big) &=& \bz(\bx) &\text{on } \,\, \Omega\times\cD \\
        \bu(\omega, \bx) &=& 0 &\text{on } \,\,  \Omega\times \partial \cD
    }\label{eq:elliptic-pde}
\end{equation}
In our numerical experiments,
we use the squared domain $\cD = (0, 1)^2$,
with the target state $\bd(\bx) = \sin(2\pi x_1) \sin(2\pi x_2)$
and the cost factor $\lambda = 10^{-8}$.
Two independent GRF realizations $\by_\ell^{(m)}$
with different mesh diameters are shown
in the four leftmost plots of \Cref{fig:log-normal-fields-and-control}.
The illustration of varying mesh resolutions
for the same sample is motivated by
telescoping sum~\eqref{eq:mlmc-estimator},
where the difference between solutions
at different levels is computed.

To find the optimal control, shown in the rightmost plots of \Cref{fig:log-normal-fields-and-control},
one common approach is to use batched stochastic gradient descent (SGD)
methods~\cite{chen2024minibatch, geiersbach2023optimization}
which, given some step size rule $t_k$, iteratively optimize the control
to find the minima of~\eqref{eq:ocp-objective}
through the scheme
\begin{equation}
    \label{eq:bsgd-iteration}
    \bz_\ell^{(k+1)} = \bz_\ell^{(k)} - t_k \nabla J(\bz_\ell^{(k)}).
\end{equation}
The gradient $\nabla J(\bz_\ell^{(k)})$ is approximated by solving
the adjoint system to~\eqref{eq:elliptic-pde} for $M_\ell$ samples
which are combined to an MC estimate.
The novelty of the work~\cite{baumgarten2025multilevel}
is that this MC estimate within the SGD method is replaced
by a MLMC estimate following~\eqref{eq:mlmc-estimator}.
Similar ideas have recently been introduced
in machine learning~\cite{weissmann2022multilevel, rowbottom2025multi}, too.
We refer again to Section~\ref{sec:numerical-experiments}
where the speedup of the proposed method is demonstrated.

All PDE systems, those for Dirichlet-Neumann averaging~\cite{kutri2024dirichlet},
as well as the state and adjoint equations,
are solved using standard linear Lagrange FE,
multigrid preconditioning, and conjugate gradient methods.

Solutions to Problem~\ref{problem:forward-uq} and Problem~\ref{problem:ocp}
can only be computed approximately,
relying on a stack of numerical methods: MC estimators, FE discretizations,
linear solvers, and optimization algorithms.
Each method introduces additional computational cost and new error components.
To quantify the total errors or both algorithms,
we use the mean squared error (MSE) for the full field estimator
\begin{equation}
    \label{eq:total-error}
    \err_{\text{MLMC}} \coloneqq
    \EE \squarelr{\norm{E^{\text{ML}}[\bu_\ell] - \EE[{\bu}]}_V^2}
    \quad \text{and} \quad
    \err_{\text{MLSGD}} \coloneqq \EE \squarelr{\bnorm{\bz_L^{(K)} - \bz^*}_{W}}
\end{equation}
measuring the distance between the intractable solution $\bz^*$ of Problem~\ref{problem:ocp},
and its approximation $\bz_L^{(K)}$ after $K$ iterations on level $L$.

Key results in~\cite{baumgarten2025budgeted, baumgarten2025multilevel}
are obtained by minimizing either error term $\err_*$
in~\eqref{eq:total-error},
subject to limited computational resources,
given through a CPU-time budget $\rP \cdot \rT_0$
and a memory budget ${\mathrm{Mem}}_0$.
This naturally leads to the following knapsack problem:

\begin{problem}[Knapsack Problem]
    \label{problem:knapsack-ocp}
    Find an optimal sequence $\bset{\tset{M_{k,\ell}}_{\ell=0}^{L_k}}_{k=0}^{K-1}$
    such that either error of~\eqref{eq:total-error} is minimized,
    while staying within a CPU-time budget $\rP \cdot \rT_{0}$ and memory ${\mathrm{Mem}}_{0}$ constraint:
    \begin{subequations}
        \label{eq:knapsack-ocp}
        \begin{align}
            \label{eq:knapsack-ocp-error}
            \hspace{-0.5cm} \min_{\,\,\,\,\,\set{\tset{M_{k,\ell}}_{\ell=0}^{L_k}}_{k=0}^{K-1}}
            \err_* \quad &\text{such that:} \\[1mm]
            \label{eq:knapsack-ocp-ct-constraints}
            \sum_{k=0}^{K_{\phantom{,}}\!\!-1}
            \sum_{\ell = 0}^{L_{k\phantom{,}}}
            \sum_{m=1}^{M_{k,\ell}} \rC^{\mathrm{CT}}_{m,\ell,k} \leq \rP \cdot \rT_{0}
            \quad &\text{and} \quad
            \rC^{\mathrm{Mem}}_{L,K-1} < {\mathrm{Mem}}_{0} \,.
        \end{align}
    \end{subequations}
    Here $\rC^{\mathrm{CT}}_{m,\ell,k}$ is the CPU-time cost of computing
    the $m$-th sample on level $\ell$ at iteration $k$,
    and $\rC^{\mathrm{Mem}}_{L,K-1}$ is the memory cost
    of storing the full field solution on level $L$ after $K-1$ iterations.
\end{problem}

\begin{figure}[t]
    \centering
    \includegraphics[width=0.8\textwidth]{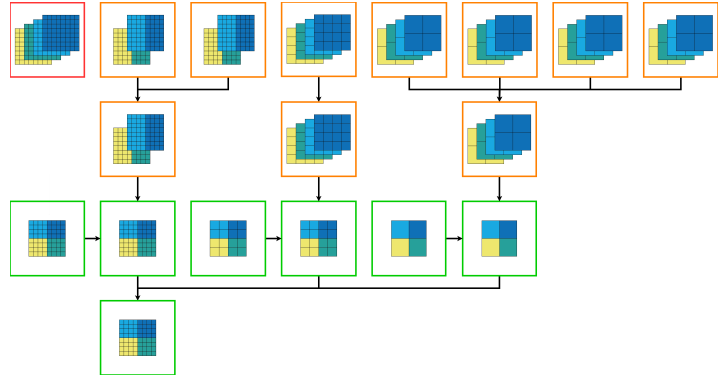}
    \caption{Multiindex FE mesh used for the MLSGD and the MLMC methods.}
    \label{fig:full-solution-update-example}
\end{figure}

We solve these problems using distributed dynamic programming techniques,
i.e., a parallel recursive algorithm which adapts in each iteration $k$
the sample sequence $\set{M_{k,\ell}}_{\ell=0}^{L_k}$ based on
previously collected data.
For implementation details, see~\cite{baumgarten2024fully, baumgarten2025budgeted, baumgarten2025multilevel};
for an intuitive overview of the multiindex data structure,
refer to Figure~\ref{fig:full-solution-update-example}.
The figure shows FE meshes at three discretization levels,
combining domain and sample parallelization.
The colors indicate the assignment of samples
and subdomains to MPI ranks.

The telescoping sum~\eqref{eq:mlmc-estimator} enables a
decreasing number of samples $\set{M_\ell}_{\ell=0}^L$
from lower to higher levels,
which enables the speedups observed in Section~\ref{sec:numerical-experiments}.
The data structure mirrors this pattern,
favoring sample parallelization on lower levels
and domain parallelization on higher levels,
ensuring that memory usage is dominated by the finest level
while minimizing inter-rank communication.

Only the domain-parallel components (green frames) remain in memory permanently;
temporary allocations are marked with orange frames,
a measure taken to reduce memory usage and to reset
the datastructures for the next iteration in the optimization algorithm.
The red frame marks a non-feasible allocation of the highest level
in a complete sample parallel data structure.

    \section{Numerical Experiments}\label{sec:numerical-experiments}

In this section, we present recent numerical results obtained
on the CPU partition of the Horeka supercomputing system,
using between $\rP=64$ and $\rP=1,024$ processing units.
The complete development, from testing, benchmarking,
and collecting the final results,
required only approximately 170,000 CPU hours
and was performed via the CI/CD pipeline described in~\cite{baumgarten2025continuous}.
This is primarily due to the results still being limited
to two spatial dimensions and the exceptional
performance of the algorithms introduced
in~\cite{baumgarten2025budgeted, baumgarten2025multilevel}.
Preliminary tests in three spatial dimensions have also
been successfully conducted with up to $\rP=16,384$ processing units,
demonstrating the feasibility of the approach in these domains as well
(we refer to the outlook for an illustration in~\Cref{fig:3d-outlook}),
but focus within this section on the two-dimensional case.

A direct consequence of using multilevel methods
to solve knapsack problems like
Problem~\ref{problem:knapsack-ocp}
is that the achievable error $\err_*$
is bounded from above by the CPU-time budget
and from below by the memory constraint.

\begin{corollary}[Lower and Upper bound of MLSGD]
    \label{cor:upper-and-lower-bound}
    The minimum~\eqref{eq:knapsack-ocp-error} is bounded
    through the imposed constraints in~\eqref{eq:knapsack-ocp-ct-constraints} as
    \begin{equation}
        \label{eq:upper-and-lower-bound}
        {\mathrm{Mem}}_{0}^{-\alpha}
        < \,\, \err_* \,\,  \lesssim \,\,
        (1 - \lambda_{\mathrm{p}}) \rT_{0}^{-\delta} + \lambda_{\mathrm{p}} (\rP \cdot \rT_{0})^{-\delta}.
    \end{equation}
    Here, $\delta = \min \bset{\tfrac12, \frac{\alpha}{2 \alpha + (\gamma - \beta)}}$
    is the convergence rate with respect to the computational resources,
    $\lambda_{\mathrm{p}} \in [0, 1]$ is the parallelizable percentage of the code,
    and $\alpha, \beta, \gamma$ are the exponents describing
    the decay of the numerical error, sampling error,
    and the growth of the computational cost.
\end{corollary}

The following numerical results estimate the exponents
$\alpha, \beta, \gamma$ and the convergence rate $\delta$.
Further, the solution to Problem~\ref{problem:knapsack-ocp} is presented,
including the optimal sample sequence $\set{M_\ell}_{\ell=0}^L$
and the distribution of computational cost across levels.

\subsection{Full Field Multilevel Monte Carlo}\label{subsec:full-field-multilevel-monte-carlo}

\begin{figure}
    \centering

    \vspace{-0.3cm}

    \includegraphics[width=\textwidth]{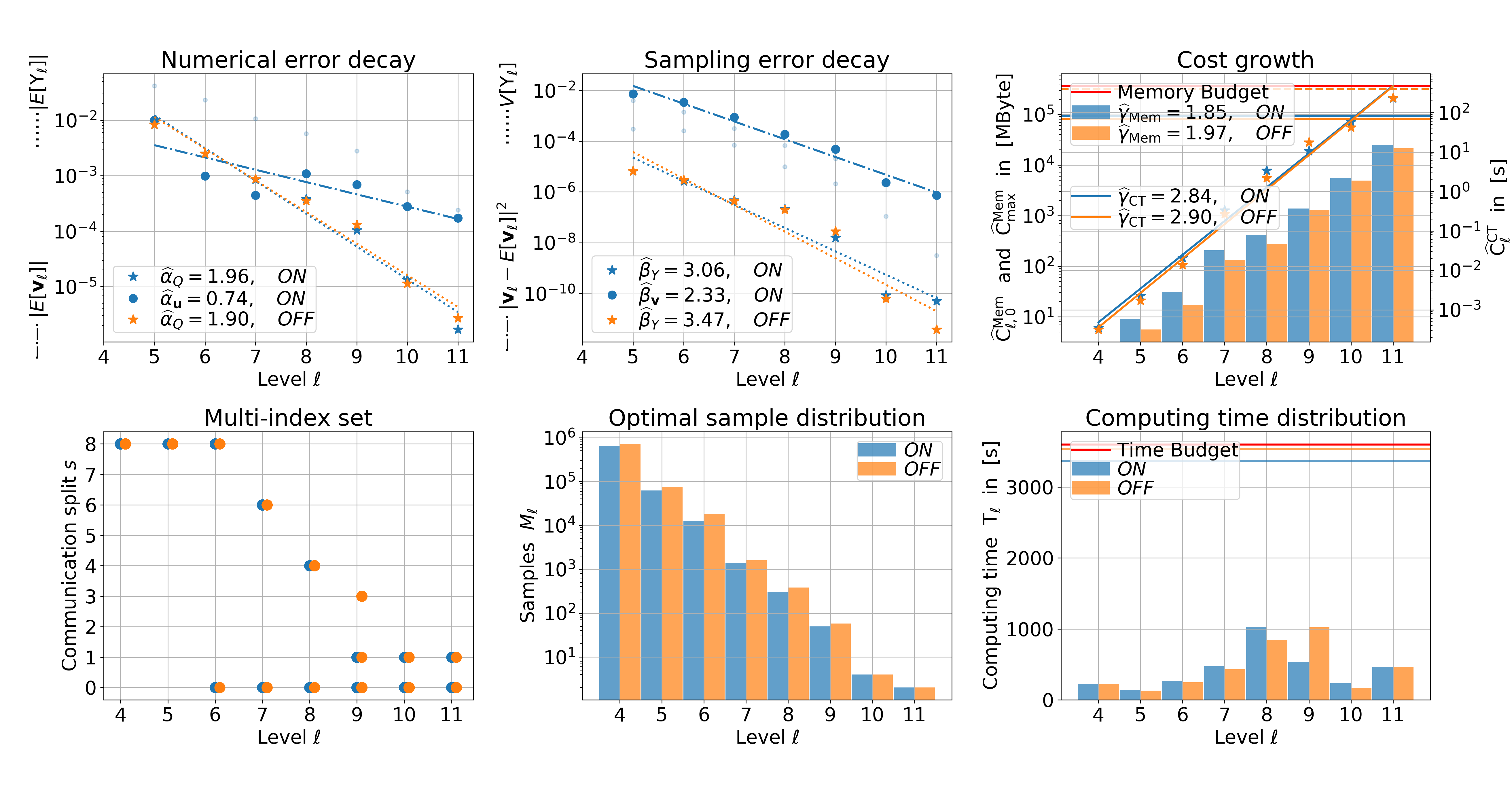}

    \vspace{-0.4cm}

    \caption{Implementations with (ON) and without (OFF) full solution update funcitonality.}
    \label{fig:computational-results-1}
\end{figure}

A central finding in~\cite{baumgarten2025budgeted} is
that estimating the full field solution,
despite theoretically weaker convergence rates $\alpha, \beta$
(see the upper left plots in~\Cref{fig:computational-results-1}),
adds no extra cost in memory footprint or CPU time
compared to estimating a scalar-valued random quantity of interest,
e.g.~$\EE[\norm{\bu}_V] \in \RR$.
This is tested at hand of the hyperbolic transport example
introduced in~\Cref{sec:numerical-experiments} on four nodes
and with $\rP = 256$ CPUs for one hour.
Although full field estimation involves substantial
data exchange between processing units,
the multiindex data structure shown
in~\Cref{fig:full-solution-update-example}
enables this communication efficiently
and without significant overhead.
This is illustrated in the upper right plot
of~\Cref{fig:computational-results-1}:
the memory usage (bar plots) and the CPU-time (solid lines)
for the full field update (highlighted in blue)
remain nearly identical to those for the scalar case (highlighted in orange).
Furthermore, the algorithm adheres to the memory
and CPU-time constraint~\eqref{eq:knapsack-ocp-ct-constraints},
visible in the upper and lower right plots of~\Cref{fig:computational-results-1}
where the total cost in memory and CPU-time (horizontal lines) stay
just below the imposed limits (horizontal red lines).
Lastly, we note the decreasing sequence of sample counts $\set{M_\ell}_{\ell=0}^L$
form the lower to the higher levels in the centered plot of the second row,
which presents the solution to Problem~\ref{problem:knapsack-ocp}.

\subsection{Multilevel Stochastic Gradient Descent}\label{subsec:multilevel-stochastic-gradient-descent}

For the MLSGD method, we present two numerical results.
The first experiment compares the MLSGD method
with the baseline batched SGD method,
where as the second experiment
demonstrates the scaling properties
of the MLSGD method.

\begin{figure}
    \includegraphics[width=1.0\textwidth]{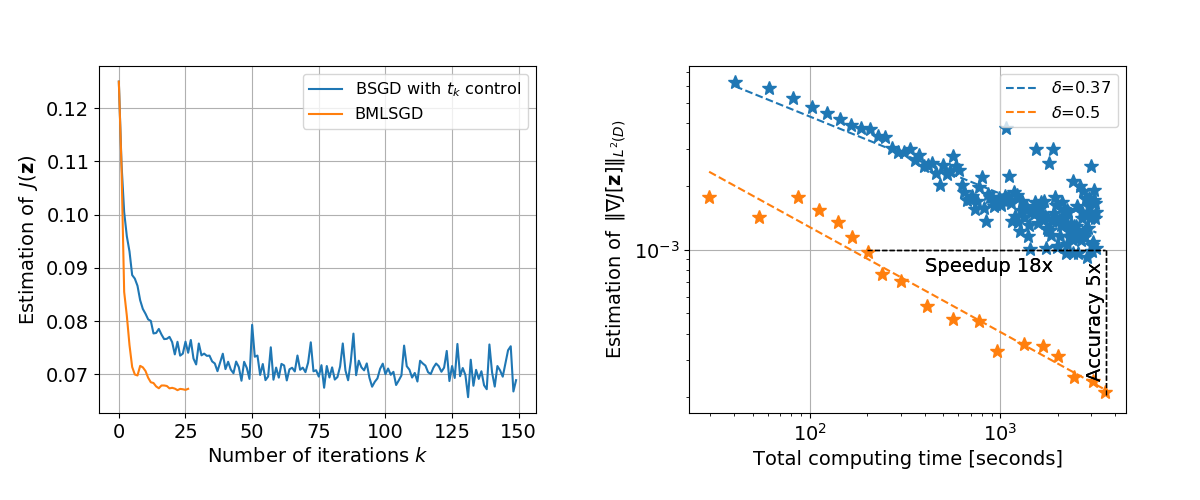}
    \vspace*{-0.5cm}
    \caption{Comparison of batched SGD and MLSGD.}
    \label{fig:metod-comparison}
\end{figure}

\subparagraph{Comparison of SGD and MLSGD}
\Cref{fig:metod-comparison} compares the MLSGD algorithm
with a standard batched SGD method on $\rP=64$ CPUs.
For SGD, we applied step size control and ran
$K = 150$ iterations of~\eqref{eq:bsgd-iteration}
using a mesh resolution of $h_\ell = 2^{-7}$
and a batch size of $M_\ell = 256$.
This configuration yielded the best SGD
performance within a one-hour computational budget.
In contrast, the MLSGD method, solving Problem~\ref{problem:knapsack-ocp}
with an optimal policy and adaptive multilevel batch sizes,
clearly outperforms SGD.
It achieves comparable accuracy 18× faster and reduces
the error by a factor of 5 at the same cost.
As shown in \Cref{fig:metod-comparison},
MLSGD also attains a better convergence rate
of $\delta \approx 0.5$, compared to $\delta \approx 0.37$ for SGD,
consistent with theoretical predictions (see~\cite{baumgarten2025multilevel}).
MLSGD also delivers improved results with fewer iterations,
lower bias, and reduced variance (see left plot of~\Cref{fig:metod-comparison}).

\begin{figure}
    \includegraphics[width=1.0\textwidth]{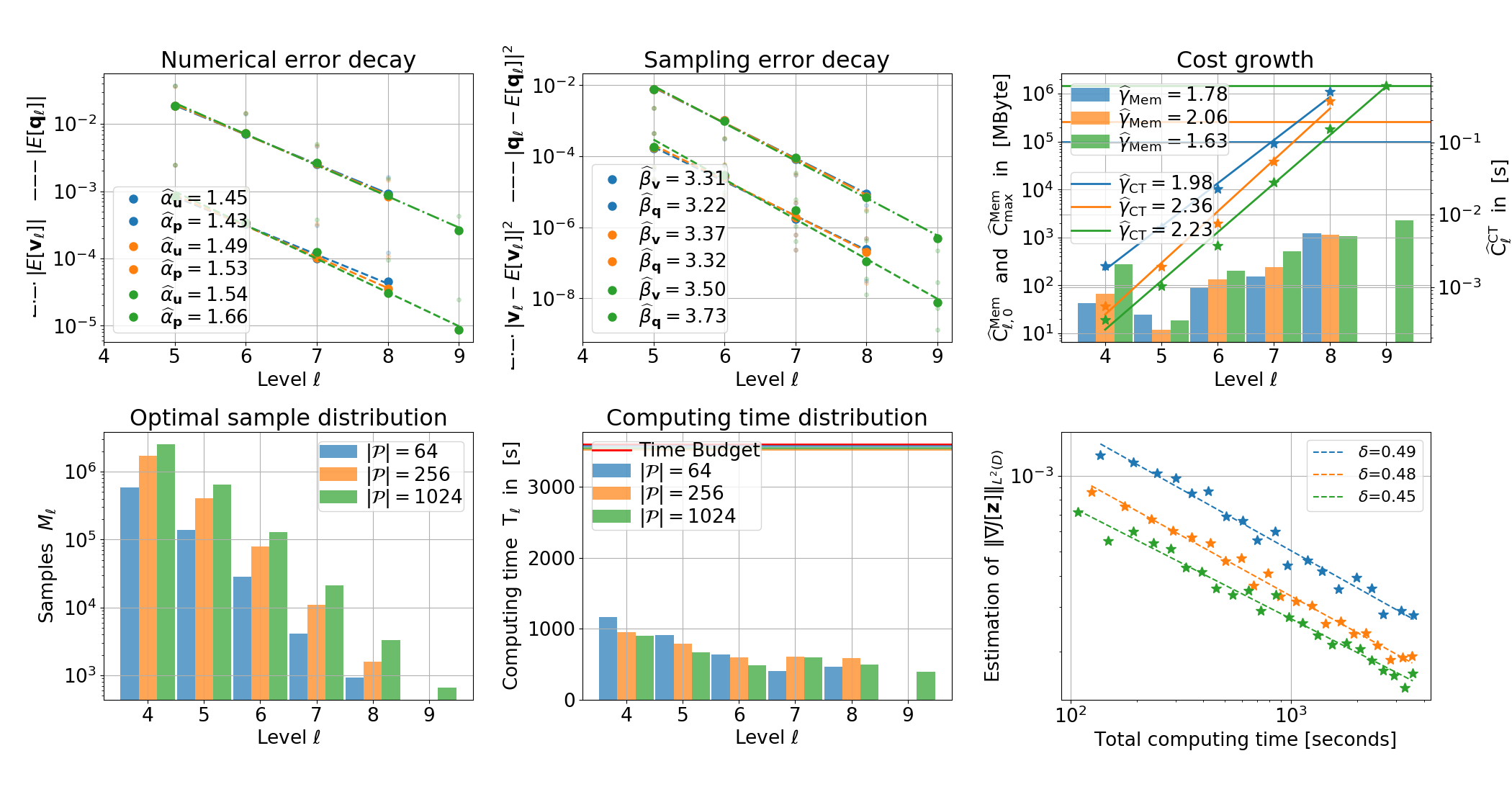}
    \vspace*{-0.5cm}
    \caption{CPU-scaling experiment with $\rP=64$, $\rP=256$ and $\rP=1024$.}
    \label{fig:nodes}
\end{figure}

\subparagraph{Node Scaling Experiment}
Finally, we demonstrate that the MLSGD method scales effectively
with increased computational resources (see~\Cref{fig:nodes}).
We run the method with the adaptive step
size rule of~\cite{koehne2024adaptivestepsizespreconditioned}
on $\rP = 64$ (1 node), $\rP = 256$ (4 nodes), and $\rP = 1,024$ (16 nodes).
The method utilizes the additional resources to compute more samples,
and, for $\rP = 1,024$, also adds an extra level
(cf.~lower left plot of~\Cref{fig:nodes},
showing the total sample count $M_\ell$ during optimization).
As illustrated in the lower right plot
and consistent with~\cite{baumgarten2024fully},
increased resources improve solution quality.
However, the smaller gap between the green ($\rP = 1,024$)
and orange ($\rP = 256$) curves,
compared to the gap between orange and blue ($\rP = 64$),
indicates diminishing parallel efficiency
This is due to the serial fraction $\lambda_{\rp}$
of the code (see~\Cref{cor:upper-and-lower-bound} and~\cite{baumgarten2024fully} for details).

    \section{Outlook}\label{sec:outlook}

\begin{figure}
    \includegraphics[trim={8cm 0 8cm 0}, clip, width=0.3\textwidth]{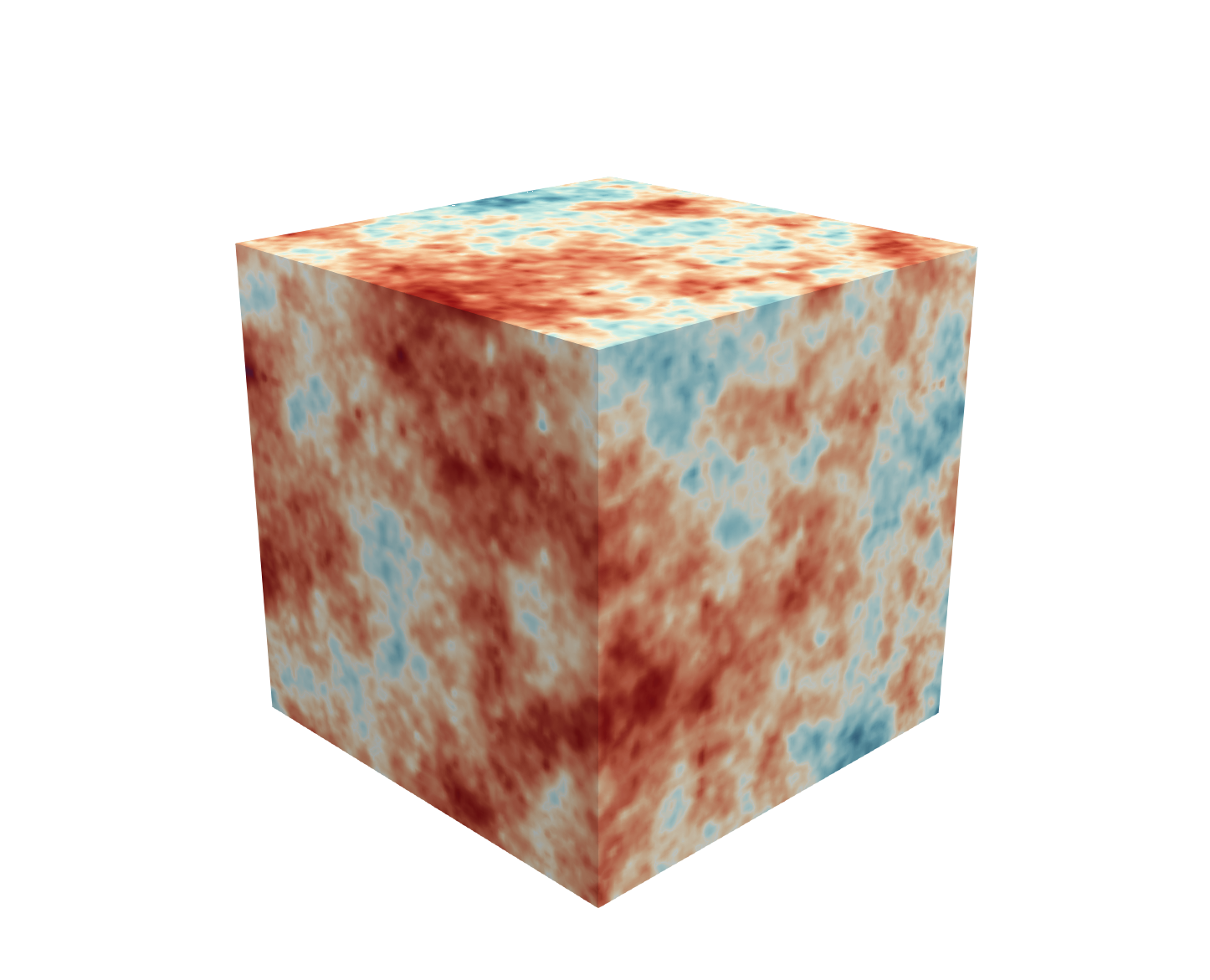}
    \includegraphics[trim={8cm 0 8cm 0}, clip, width=0.3\textwidth]{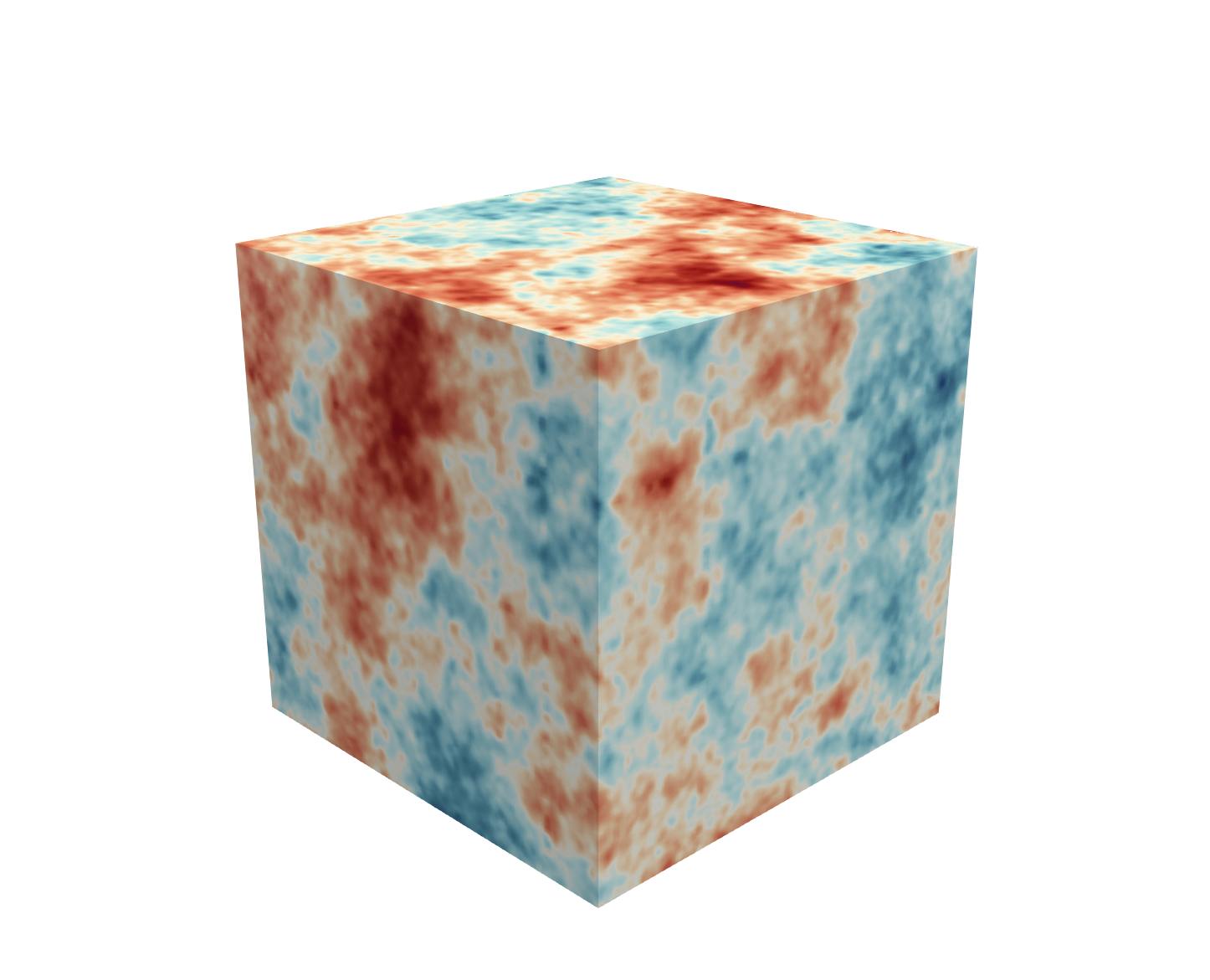}
    \includegraphics[trim={8cm 0 8cm 0}, clip, width=0.3\textwidth]{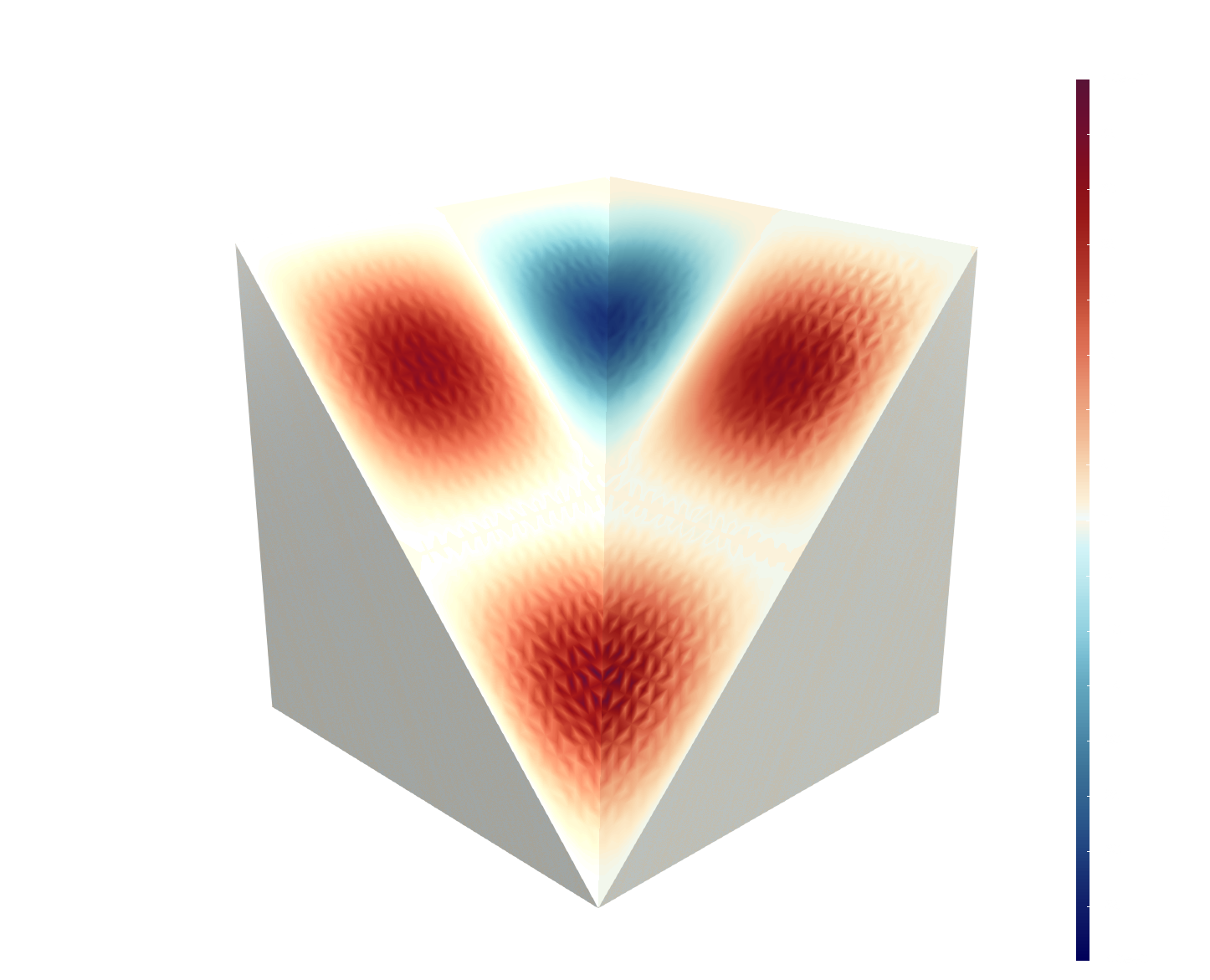}
    \vspace*{-0.5cm}
    \caption{Outlook on three-dimensional random field sampling and optimal control.}
    \label{fig:3d-outlook}
\end{figure}

The orange line in~\Cref{fig:metod-comparison} establishes
a new baseline for evaluating future developments of the MLSGD method.
Any improvements to the M++ software~\cite{baumgarten2021parallel}
can be measured against this baseline to assess their impact
in a completely automated manner following~\cite{baumgarten2025continuous}.
Software engineering and performance optimization
reduce the constant factor in convergence plots,
as shown in~\Cref{cor:upper-and-lower-bound}.
To further improve performance,
we plan to extend the implementation to GPUs
using Ginkgo~\cite{anzt2022ginkgo},
a linear algebra framework for high-performance computing.
Initial steps include designing simple test cases
and transferring data from our algebraic structures to Ginkgo’s format.
Although the current implementation is still rudimentary,
we expect that GPU parallelism will significantly
accelerate the method and further
lower the orange line in~\Cref{fig:metod-comparison}.
In the second half of the project year,
we will also focus on scaling the algorithm
to three-dimensional applications, as illustrated in~\Cref{fig:3d-outlook}.

    \begin{acknowledgement}
    The authors gratefully acknowledge the computing time provided
    on the high-performance computer HoreKa by the National High-Performance Computing Center at KIT (NHR@KIT).
    This center is jointly supported by the Federal Ministry of Education
    and Research and the Ministry of Science, Research
    and the Arts of Baden-Württemberg,
    as part of the National High-Performance Computing (NHR)
    joint funding program.
    HoreKa is partly funded by the German Research Foundation (DFG).
\end{acknowledgement}

    \bibliographystyle{ieeetr}

    \bibliography{bibliography}

\end{document}